\documentclass[journal ]{new-aiaa}
\usepackage[utf8]{inputenc}
\usepackage{textcomp}

\usepackage{graphicx}
\usepackage{subcaption}
\usepackage{amsmath}
\usepackage{commath}
\usepackage[version=4]{mhchem}
\usepackage{siunitx}
\usepackage{longtable,tabularx}
\usepackage{cancel}
\usepackage{mathtools}

\usepackage{amsthm}

\newtheorem{remark}{Remark}
\newtheorem{assumption}{Assumption}

\newcommand{\comment}[1]{}

\newcommand\blfootnote[1]{%
  \begingroup
  \renewcommand\thefootnote{}\footnote{#1}%
  \addtocounter{footnote}{-1}%
  \endgroup
}

\title{
Convex Bound of Nonlinear Dynamical Errors for Covariance Steering
\blfootnote{An earlier version of this paper was presented as paper 25-504 at the 2025 AAS/AIAA Astrodynamics Specialist Conference.}
}

\author{Daniel C. Qi\footnote{Ph.D. Student, School of Aeronautics and Astronautics. (Corresponding Author: qi85@purdue.edu)} and Kenshiro Oguri\footnote{Assistant Professor, School of Aeronautics and Astronautics. Senior Member AIAA.}}
\affil{Purdue University, West Lafayette, Indiana, 47907}

\begin{document}

\maketitle

\section{Introduction}\label{sec: intro}
A common approach for controller design in nonlinear systems is to linearize the dynamics about a reference trajectory and apply linear control techniques within the vicinity of the reference \cite{Nguyen-linear-control-in-NL}. Naturally, the omitted nonlinear terms from the linearization process will result in some inconsistencies between the expected and actual responses of the system. The degree of inconsistency is important as there is a point from the reference where the linearized dynamics no longer accurately represent the actual nonlinear dynamics, and thus any linear controller operating in this region could potentially behave unpredictably. As a result, ensuring the validity of the designed linear controllers in nonlinear environments is an important aspect of mission assurance. 

Covariance control, or covariance steering, is a stochastic optimal control problem focusing on creating a control policy for a Gaussian state distribution \cite{Hotz-Cov-Steering-OG, Okamoto-Cov-Steering}. A key assumption underlying covariance steering is that the uncertainty remains Gaussian under linear dynamics. This follows from the fact that Gaussian distributions are preserved under affine transformations. Moreover, since a Gaussian distribution is fully characterized by its mean and covariance, the stochastic optimal control problem can be reformulated as a mean and covariance control problem. For nonlinear systems, linearization techniques can be employed to extend this framework, and this has been successfully demonstrated in problems such as trajectory optimization \cite{Ridderhof-SCP, Naoya-Sequential-Cov-Steering, Jerry-CovSteering, Benedikter-CovSteering}, orbit stationkeeping \cite{Divija-stationkeeping, Oguri-Chance-Paper}, and proximity operations \cite{Geller-RPO, Philip-prox-op}. However, in highly nonlinear systems like cislunar space, originally Gaussian uncertainties often became non-Gaussian due to the nonlinear dynamics \cite{erin-cislunar-OD, John-GMM, Sharad-UQ-CR3BP}. Thus, linear covariance controllers have shown discrepancies between the predicted covariance and the actual distribution observed in nonlinear simulations \cite{Naoya-Sequential-Cov-Steering, Jerry-CovSteering}. 

This leads to a fundamental challenge in applying covariance steering to nonlinear systems. Unlike deterministic trajectory optimization, where sequential methods iteratively refine a nominal trajectory through successive linearizations, covariance steering must also ensure that the linearization accurately captures the evolution of both the mean and the covariance. While sequential techniques can be effectively applied to the mean steering component, the covariance evolution is intrinsically coupled to the linearized dynamics, which are defined about the mean trajectory \cite{Naoya-Sequential-Cov-Steering, Benedikter-CovSteering}. This coupling prevents independent refinement of the covariance and precludes a straightforward extension of sequential methods to the covariance component.

One way to address this issue is to limit the distribution’s exposure to regions of strong nonlinearity. For covariance steering, this equates to reducing the covariance so that it remains well-approximated by a linear transformation. By keeping the distribution concentrated, a valid linear approximation can better capture its behavior, making linear covariance control methods more effective even in inherently nonlinear systems. But to avoid highly nonlinear regions, they must first be identified. This raises the idea of quantifying the nonlinearity of a region. Determining Measures of Nonlinearities (MoNs) has been explored in the past \cite{Beale-MoN, junkins-MoN}. A drawback of many MoNs is the necessity for empirical sampling which brings a computational cost that is impractical for any controller optimization scheme. Semi-analytical methods such as the Tensor Eigenpair Measure of Nonlinearity \cite{Jenson-MeasuresOfNonlinearity} do exist, and it is the primary inspiration for this paper's approach. However, these MoNs are usually only investigated from a dynamics standpoint. In the context of controls, previous MoN-based controllers use MoNs to determine when the system should be re-linearized about a new reference point \cite{Nguyen-linear-control-in-NL, Chen-MoN-MPC}. However, these frameworks primarily serve as analysis tools for detecting the need for re-linearization, rather than integrating this mechanism directly into the trajectory planning process for problems such as covariance steering.

Developing a covariance steering framework that accounts for nonlinearities can extend the controller's time horizon and enhance its ability to perform in nonlinear systems. Most existing MoNs are not formulated in a way that enables directly constraining or minimizing nonlinearities within an optimization framework. In other words, they primarily quantify local, instantaneous nonlinearity at a given point, but do not provide a mechanism for regulating it over a time horizon. This paper presents a formulation of an MoN that enables optimization-based covariance controllers to directly minimize dynamical nonlinearities. The approach involves calculating higher-order dynamical terms and then bounding the contributions of their nonlinear effects. This bound is then used as an objective function for a control optimization problem. A crucial aspect is that this formulation is convex and readily solvable with standard convex optimization solvers. As a result, the controller minimizes its exposure to regions of high nonlinearity, thereby enabling covariance steering control policies to be extended to more nonlinear regimes. The effectiveness of this method has been successfully demonstrated in spacecraft trajectory optimization under uncertainty \cite{Frassinella-min-NL}. However, the present work constitutes the original and comprehensive presentation of the approach.

\section{Background}\label{sec: background}
\subsection{Tensor Operations}
This section details the notation used for tensor operations and is mainly borrowed from previous works \cite{Jenson-MeasuresOfNonlinearity}. Consider an arbitrary real $\ell$-th order tensor $\mathcal{T}^{(\ell)}$ with dimensions $n$ along each of its indices. This leads to a tensor with $n^\ell$ elements with each element denoted by
\begin{equation}
    \mathcal{T}_{i_1 i_2 \ldots i_\ell} \quad \text{for} \quad 1 \leq i_1,i_2,\ldots,i_\ell \leq n
\end{equation}

Let $\boldsymbol{v}$ be an arbitrary real vector with dimensions $n$ whose $i$-th element is expressed as $v_i$. Multiplying $\mathcal{T}^{(\ell)}$ with $\boldsymbol{v}$ in $\ell-1$ modes result in a $n$ dimensional vector. This product is denoted by $\mathcal{T}^{(\ell)} \cdot \boldsymbol{v}^{\ell-1}$, and the $j$-th element of this new vector can be expressed using a summation convention
\begin{equation}
    \left( \mathcal{T}^{(\ell)} \cdot \boldsymbol{v}^{\ell-1} \right)_j = 
    \sum\limits_{i_2 \ldots i_\ell} \mathcal{T}_{j i_2 \ldots i_\ell} \underbrace{v_{i_2} \ldots v_{i_\ell}}_{\ell-1} 
\end{equation}
where $\sum_{i_2 \ldots i_\ell}$ denotes a series of nested summations $\sum^n_{i_2 = 1} \ldots \sum^n_{i_\ell = 1}$. 

Note that multiplying an $n \times n$ matrix $A$ with a vector is equivalent to a 2nd-order tensor multiplied by a vector. The $j$-th element of the Matrix-vector product is simply 
\begin{equation}\label{eq: tensor, matrix times vector}
    \left(A \boldsymbol{v} \right)_j = 
    \sum\limits_{i_2}^{n} A_{j i_2} v_{i_2} 
\end{equation}

\subsection{State Transition Matrices and Tensors}
The dynamics around a reference trajectory can be approximated with a series expansion using the state transition matrix (STM) and the state transition tensors (STT) \cite{Turner-STT,Pellegrini-STT,Jenson-MeasuresOfNonlinearity}. Given a state $\boldsymbol{\mathcal{X}}$ and some reference $\boldsymbol{\mathcal{X}}^*$, its dynamics can be written in  the following form:
\begin{equation}
\begin{aligned}
    \dot{\boldsymbol{\mathcal{X}}} = f(t,\boldsymbol{\mathcal{X}}) &\quad \rightarrow \quad &\dot{\boldsymbol{\mathcal{X}}} = \dot{\boldsymbol{\mathcal{X}}}^* + \delta\dot{\boldsymbol{\mathcal{X}}}  = f(t,\boldsymbol{\mathcal{X}}^* + \delta\boldsymbol{\mathcal{X}})
\end{aligned}
\end{equation}
where $\delta\boldsymbol{\mathcal{X}} = \boldsymbol{\mathcal{X}} - \boldsymbol{\mathcal{X}}^*$. Using a series expansion on the dynamics about $\boldsymbol{\mathcal{X}}^*$, integrating from $t_k$ to $t_{k+1}$, the state deviation from the reference at $t_{k+1}$ can be written as a function of the state deviation at $t_{k}$ using the STM and STTs.
\begin{equation}\label{eq: dynamics series expansion}
\begin{aligned}
	\delta\boldsymbol{\mathcal{X}}_{k+1}
    &= \Phi(t_{k+1},t_k) \delta\boldsymbol{\mathcal{X}}_{k} + \sum\limits^{\infty}_{m = 2} \frac{1}{m!} \Phi^{(m)}(t_{k+1},t_k) \cdot \delta\boldsymbol{\mathcal{X}}^m_{k}
\end{aligned}
\end{equation}
where $\delta\boldsymbol{\mathcal{X}}_{k}$ is the state deviation about the reference at time $t_k$, $\Phi(t_{k+1},t_k)$ is the first-order STM from $t_k$ to $t_{k+1}$, and $\Phi^{(m)}(t_{k+1},t_k)$ is the $m$-th order STM or STT. Note the different indexing between STM and tensors (i.e., $\ell = m +1$): since $\Phi$ is a matrix, it is a second-order tensor. Then, $\Phi^{(2)}$ is a third-order tensor, etc. 

\section{Impulsive Linear Covariance Steering}
\subsection{Problem Setup for Impulsive Linear Covariance Steering}\label{sec: cov steering}
This paper considers the finite-horizon, discrete-time covariance steering problem. The objective of covariance steering is to determine an optimal control policy that simultaneously regulates the mean and covariance of a state distribution. Such problems can be formulated and solved using convex optimization techniques \cite{Okamoto-Cov-Steering, Liu-CS}, with numerous extensions developed to address different applications. In particular, one important variant is the output-feedback formulation, where the objective is to control the true state distribution when control actions are based only on noisy state estimates \cite{Ridderhof-CS, Pilipovsky-CS-Output-Feedback, Divija-stationkeeping, Naoya-Sequential-Cov-Steering}. This approach has been proposed for spacecraft autonomy due to its more realistic integration of estimation and control, as discussed by Oguri~\cite{Oguri-Chance-Paper}. There exist several formulations of the output-feedback covariance steering problem, depending on how the state statistics are computed. This paper adopts the formulation from Oguri~\cite{Oguri-Chance-Paper}, which employs a block Cholesky parameterization for linear covariance steering. However, the proposed approach for minimizing MoN is not restricted to this formulation and can be readily extended to other parameterizations, as it does not depend on their specific structure.

\begin{assumption}
An additional assumption in this work is the restriction to affine impulsive control actions. This is motivated by the fact that the proposed MoN formulation captures only dynamical nonlinearities, and does not account for nonlinearities arising from non-affine control mappings or control interpolation schemes such as zero-order or first-order hold methods \cite{Hofmann-ZOH}. Under impulsive control, the state is updated instantaneously at each control step, so that linearization over each interval introduces nonlinear error only from the system dynamics.
\end{assumption}

Let the time horizon be discretized into $N$ nodes, forming $N-1$ segments, with $k = 0,1, \ldots, N-1$ denoting the time index. Let $\boldsymbol{x}_{k} \in \mathbb{R}^{n_x}$ a multivariate random variable denoting the \emph{true} state, and $\boldsymbol{u}_{k} \in \mathbb{R}^{n_u}$ be the impuslive control action. Consider the following stochastic systems under impulsive control: 
\begin{equation}
\begin{aligned}
    \boldsymbol{x}_{k+1} &= A_k \boldsymbol{x}^{+}_k +\boldsymbol{c}_k + G_k \boldsymbol{w}_k 
    &\qquad &\text{(Dynamics Update)}\\
    \boldsymbol{y}_{k+1} &= C_{k+1} \boldsymbol{x}_{k+1} + D_{k+1}\boldsymbol{w}_{\text{obs},k+1} + \boldsymbol{c}_{\text{obs},k+1} 
    &\qquad &\text{(Measurement Model)}\\
    \boldsymbol{x}^{+}_{k+1} &= \boldsymbol{x}_{k+1} + B\boldsymbol{u}_{k+1} 
    &\qquad &\text{(Impuslive Control Input)}
\end{aligned}
\end{equation}
where $\boldsymbol{x}_{k}^+ \in \mathbb{R}^{n_x}$ is the state vector after the control action and $B\in \mathbb{R}^{n_x \times n_u}$ is the affine control mapping matrix. Note the distinction between $\boldsymbol{x}_{k}$ and $\boldsymbol{x}^{+}_{k}$ to denote pre- and post- maneuver states at a given time $t_k$. Additionally, $A_k\in \mathbb{R}^{n_x\times n_x}$ and $\boldsymbol{c}_k\in \mathbb{R}^{n_x}$ are \emph{linearized} state mappings from the \emph{nonlinear} dynamics, $G_k\in \mathbb{R}^{n_x\times n_x}$ captures any miscellaneous stochastic events that affect the state, and $\boldsymbol{w}_k\in \mathbb{R}^{n_x}$ is a standard Gaussian white noise \cite{Oguri-Chance-Paper}. For the measurement model, $\boldsymbol{y}_{k} \in \mathbb{R}^{n_y}$ represents the measurement taken at $t_k$, and $C_k \in \mathbb{R}^{n_y \times n_x}$ and $\boldsymbol{c}_{\text{obs},k} \in \mathbb{R}^{n_y}$ are the linearized measurement mappings, and $D_k \in \mathbb{R}^{n_y \times n_y}$ is the measurement error matrix with corresponding standard Gaussian white noise $\boldsymbol{w}_{\text{obs},k} \in \mathbb{R}^{n_y}$ \cite{Oguri-Chance-Paper}. Now let $\hat{\boldsymbol{x}}^-_{k} \in \mathbb{R}^{n_x}$ denote the a priori state estimate, $\hat{\boldsymbol{x}}_{k} \in \mathbb{R}^{n_x}$ denote the a posteriori state estimate before impulsive control, and $\hat{\boldsymbol{x}}^+_{k} \in \mathbb{R}^{n_x}$ denote the a posteriori state estimate after impulsive control. The estimated state has the following updates following standard Kalman filtering:
\begin{equation}\label{state updates}
    \begin{aligned}
    \hat{\boldsymbol{x}}^{-}_{k+1} &= A_{k}\hat{\boldsymbol{x}}^+_{k} + \boldsymbol{c}_k
    &\qquad &\text{(Filter Time Update)}\\
    \hat{\boldsymbol{x}}_{k+1} &= \hat{\boldsymbol{x}}^{-}_{k+1} + L_{k+1} \tilde{\boldsymbol{y}}^{-}_{k+1}
    &\qquad &\text{(Filter Measurement Update)} \\
   \hat{\boldsymbol{x}}^+_{k+1} &= \hat{\boldsymbol{x}}_{k+1} + B \boldsymbol{u}_{k+1}
    &\qquad &\text{(Impulsive Control Input)} 
    \end{aligned}
\end{equation}
where $\tilde{\boldsymbol{y}}^{-}_{k}$ is the innovation process from Kalman filtering. The Kalman gain is calculated with the following:
\begin{equation}
    L_k = \tilde{P}^{-}_k C_k^\top (C_k \tilde{P}^{-}_k C_k^\top +D_k D_k^\top)^{-1}
\end{equation}
where the time and measurement update covariances about the reference can be calculated. 
\begin{equation} \label{eq: Ptilde}
\begin{aligned}
    \tilde{P}^{-}_k = A_{k-1} \tilde{P}_{k-1}A_{k-1}^\top + G_{k-1}G_{k-1}^\top, &\qquad &
    \tilde{P}_k = (I-L_k C_k)\tilde{P}^{-}_k(I-L_k C_k)^\top + L_k D_k D_k^\top L_k^\top
\end{aligned}
\end{equation}
Note that the Kalman gain can be calculated \emph{before} the time and measurement updates as long as the measurement error covariance is known \cite{Oguri-Chance-Paper, Ridderhof-CS, Pilipovsky-CS-Output-Feedback, Naoya-Sequential-Cov-Steering}. 

\subsection{Block Cholesky Formulation for Impulsive Linear Covariance Steering}\label{sec: impulsive linear control}
The block Cholesky procedure follows that of Oguri~\cite{Oguri-Chance-Paper}, but is applied to post-maneuver state statistics. By combining the time, measurement, and control updates in Eq.~\eqref{state updates}, a single equation for the state estimate is obtained.
\begin{equation}
    \begin{aligned}
    \hat{\boldsymbol{x}}^+_{k+1} &= A_{k}\hat{\boldsymbol{x}}^+_{k} + B\boldsymbol{u}_{k+1} + \boldsymbol{c}_k  + L_{k+1} \tilde{\boldsymbol{y}}^{-}_{k+1}
    \end{aligned}
\end{equation}

The state and control throughout a time horizon can then be expressed in a block-matrix form.
\begin{equation}
\begin{aligned}
    \begin{bmatrix}
    \hat{\boldsymbol{x}}^+_0 \\
    \hat{\boldsymbol{x}}^+_1 \\
    \hat{\boldsymbol{x}}^+_2 \\
    \vdots
    \end{bmatrix}
    &=
    \begin{bmatrix}
    I_{n_x} \\
    A_0 \\
    A_1 A_0 \\
    \vdots
    \end{bmatrix}
    \hat{\boldsymbol{x}}^{-}_0
    + 
    \begin{bmatrix}
    B & 0 & 0 & \\
    A_0 B & B & 0 & \cdots\\
    A_1 A_0 B & A_1 B & B & \\
    & \vdots & & \ddots
    \end{bmatrix}
    \begin{bmatrix}
    \boldsymbol{u}_0 \\
    \boldsymbol{u}_1 \\
    \boldsymbol{u}_2 \\
    \vdots
    \end{bmatrix} 
    + 
    \begin{bmatrix}
    0 & 0 & 0 & \\
    I_{n_x} & 0 & 0 & \cdots\\
    A_1  & I_{n_x} & 0 & \\
    & \vdots & & \ddots
    \end{bmatrix}
    \begin{bmatrix}
    \boldsymbol{c}_0 \\
    \boldsymbol{c}_1 \\
    \boldsymbol{c}_2 \\
    \vdots
    \end{bmatrix} 
    +
    \begin{bmatrix}
    L_0 & 0 & 0 & \\
    A_0 L_0 & L_1 & 0 & \cdots\\
    A_1 A_0 L_0 & A_1 L_1 & L_2 & \\
    & \vdots & & \ddots
    \end{bmatrix}
    \begin{bmatrix}
    \tilde{\boldsymbol{y}}^{-}_0 \\
    \tilde{\boldsymbol{y}}^{-}_1 \\
    \tilde{\boldsymbol{y}}^{-}_2 \\
    \vdots
    \end{bmatrix}
\end{aligned}
\end{equation}
where $I_{n}$ is an $n \times n $ identity matrix. This is expressed in a more compact form:
\begin{equation} \label{eq: dynamics block form}
    \hat{\mathbf{X}}^+ = \mathbf{A}\hat{\boldsymbol{x}}^{-}_0 + 
    \mathbf{B}^+\mathbf{U} + \mathbf{C}+
    \mathbf{L}\mathbf{Y}
\end{equation}

The formulation from Oguri~\cite{Oguri-Chance-Paper} defines the linear feedback control policy as
\begin{equation}
\begin{aligned}
   \boldsymbol{u}_k = \bar{\boldsymbol{u}}_k + K_k \boldsymbol{z}_{k}, &\qquad &
   \boldsymbol{z}_{k+1} = A_k \boldsymbol{z}_{k} + L_{k+1}\tilde{\boldsymbol{y}}^{-}_{k+1}  
\end{aligned}
\end{equation}
where $\boldsymbol{z}_0 = \hat{\boldsymbol{x}}_0 - \bar{\boldsymbol{x}}_0$. This can also be expressed in block-matrix and a compact form.

\begin{equation}\label{eq: control policy block form}
    \mathbf{U} = \bar{\mathbf{U}} + \mathbf{K}\mathbf{Z}, \qquad
    \mathbf{Z} = \mathbf{A}(\hat{\boldsymbol{x}}^{-}_0 - \bar{\boldsymbol{x}}_0) + \mathbf{L}\mathbf{Y}
\end{equation}

Eq.~\eqref{eq: dynamics block form} and \eqref{eq: control policy block form} can then be combined into one expression.
\begin{equation}
    \hat{\mathbf{X}}^+ = \mathbf{A}\hat{\boldsymbol{x}}^{-}_0 + 
    \mathbf{B}^+\bar{\mathbf{U}} + \mathbf{B}^+\mathbf{K}\mathbf{Z}+\mathbf{C}+
    \mathbf{L}\mathbf{Y} 
\end{equation}

Since $\mathbb{E}[\mathbf{Y}] = \vec{0}$, $\mathbb{E}[\mathbf{Z}] = \vec{0}$, and  $\mathbb{E}[\mathbf{U}] = \bar{\mathbf{U}}$, the expected value and covariance of the estimated state are
\begin{subequations} 
\begin{align}
    \bar{\mathbf{X}}^+ &\triangleq \mathbb{E}[\hat{\mathbf{X}}^+] = \mathbf{A}\bar{\boldsymbol{x}}_0 + \mathbf{B}^+\bar{\mathbf{U}} + \mathbf{C}  \label{eq: cov steering mean}\\
    \hat{\mathbf{P}}^+ &\triangleq \text{Cov}[\hat{\mathbf{X}}^+] = (I + \mathbf{B}^+\mathbf{K})\mathbf{S}(I + \mathbf{B}^+\mathbf{K})^\top
\end{align}
\end{subequations}
where $\mathbf{S} \triangleq \text{Cov}[\mathbf{A}(\hat{\boldsymbol{x}}^{-}_0 - \bar{\boldsymbol{x}}_0) + \mathbf{L}\mathbf{Y}]$. Noting that the initial states and measurements are independent, 
\begin{equation}
    \begin{aligned}
\mathbf{S} &= \text{Cov}[\mathbf{A}(\hat{\boldsymbol{x}}^{-}_0 - \bar{\boldsymbol{x}}_0)]+\text{Cov}[ \mathbf{L}\mathbf{Y}] =\mathbf{A} \hat{P}_{0^-}\mathbf{A}^\top + \mathbf{L}\mathbf{P}_Y\mathbf{L}^\top
    \end{aligned}
\end{equation}
where $\mathbf{P}_Y = \text{Cov}[\mathbf{Y}]=\text{blkdiag}\left( P_{\tilde{\boldsymbol{y}}^{-}_0}, P_{\tilde{\boldsymbol{y}}^{-}_1},\ldots \right)$ and $P_{\tilde{\boldsymbol{y}}^{-}_k} = C_k \tilde{P}^-_k C_k^\top + D_k D_k^\top$. Applying the lower-triangular decomposition results in an affine function in $\mathbf{K}$.
\begin{equation} \label{eq: Phat}
    \begin{aligned}
(\hat{\mathbf{P}}^+ )^{1/2} = (I +\mathbf{B}^+\mathbf{K})\mathbf{S}^{1/2}
    \end{aligned}
\end{equation}

Previous works \cite{Ridderhof-CS, Pilipovsky-CS-Output-Feedback, Naoya-Sequential-Cov-Steering} have shown that adding the covariances from Eq.\eqref{eq: Ptilde} and \eqref{eq: Phat} yields the covariance of the true state under feedback control $P^+_k = \hat{P}^+_k + \tilde{P}_k$. The lower-triangular covariance for the true state covariance is then
\begin{equation} \label{eq: cov steering cov}
\begin{aligned} 
    (P_k^+)^{1/2} = \left[(\hat{P}^+_k)^{1/2} \quad \tilde{P}_k^{1/2}\right]
\end{aligned}
\end{equation}
where $(\hat{P}^+_k)^{1/2} = E_{x_k}(\hat{\mathbf{P}}^+)^{1/2}$ and the extraction function $E_{x_k}$ defined as $\boldsymbol{x}^+_k = E_{x_k}\mathbf{X}^+$. The final convex optimization problem for the block Cholesky covariance steering becomes
\begin{equation}
\begin{aligned}
\min_{\bar{\boldsymbol{X}},\hat{\boldsymbol{P}}^{1/2},\bar{\boldsymbol{U}},\boldsymbol{K}} & 
J(\bar{\boldsymbol{X}},\hat{\boldsymbol{P}}^{1/2},\bar{\boldsymbol{U}},\boldsymbol{K})
&\text{(Convex Objective Function)}
\\
\text{s.t. }
& \text{Eq}.~\eqref{eq: cov steering mean}&\text{(Mean Dynamics)}\\
& \text{Eq}.~\eqref{eq: Phat},\eqref{eq: cov steering cov}&\text{(Covariance Dynamics)}\\
& h_{\text{eq}}(\bar{\boldsymbol{X}},\hat{\boldsymbol{P}}^{1/2},\bar{\boldsymbol{U}},\boldsymbol{K}) = 0&\text{(Affine Equality Constraint)}\\
& h_{\text{ineq}}(\bar{\boldsymbol{X}},\hat{\boldsymbol{P}}^{1/2},\bar{\boldsymbol{U}},\boldsymbol{K}) \leq 0&\text{(Convex Inequality Constraint)}\\
\end{aligned}
\end{equation}
where $h_{\text{eq}(\cdot)}$ and $h_{\text{ineq}(\cdot)}$ are any additional affine equality and convex inequality constraints imposed on the problem.

\subsection{Chance Constraint and Quantile}
Within the framework of convex optimization, convex constraints can be applied to both the controller and the state. For covariance controllers, chance constraints can be implemented. The relevant chance constraints for this paper are specified below, with more types of constraints discussed by Oguri~\cite{Oguri-Chance-Paper}. Given that $\boldsymbol{x} \sim \mathcal{N}(\bar{\boldsymbol{x}}, P) \in \mathbb{R}^n$, the probabilistic norm constraint can be upper-bounded by an expression given its mean and covariance:
\begin{equation}\label{eq: chance constraints}
\mathbb{P}\left(\norm{\boldsymbol{x} - \boldsymbol{a}}_2\leq b\right) \geq 1 - \epsilon
\qquad
\rightarrow
\qquad
\norm{\bar{\boldsymbol{x}} - \boldsymbol{a}}_2 + \sqrt{m_{\chi^2}(\epsilon,n)}\norm{P^{1/2}}_2 \leq b
\end{equation}
where $P = (P^{1/2})(P^{1/2})^\top$. The function $m_{\chi^2}(\epsilon,n)$ is the inverse of a chi-square's cumulative distribution function with $n$ degrees of freedom evaluated at the probability value of $1-\epsilon$, and $\epsilon$ is a user-defined value that determines a probabilistic ``risk'' factor. Since Gaussian distributions are unbounded, $0< \epsilon \ll 1$ determines the amount of distribution that will satisfy the constraint. A smaller $\epsilon$ increases the probability of success but may lead to solution infeasibility, as the quantile used to classify an event as safe becomes too large.

From Eq.~\eqref{eq: chance constraints}, it can be seen that the constraint can lead to an upper bound on the quantile function. Given that $\boldsymbol{x} \sim \mathcal{N}(\bar{\boldsymbol{x}}, P) \in \mathbb{R}^n$, the quantile for the $l^2$norm of $\boldsymbol{x} - \boldsymbol{a}$ evaluated at probability $1-\epsilon$ is upper bounded by $\text{Quant}_{\text{ub}}$:
\begin{equation}\label{eq: chance constraints quantile}
\text{Quant}(1-\epsilon) \leq \text{Quant}_{\text{ub}}(1-\epsilon) = \norm{\bar{\boldsymbol{x}} - \boldsymbol{a}}_2 + \sqrt{m_{\chi^2}(\epsilon,n)}\norm{P^{1/2}}_2
\end{equation}
It is important to note that Eq.~\eqref{eq: chance constraints} and \eqref{eq: chance constraints quantile} work under the $l^2$norm of the Gaussian distributions.

\section{Problem Statement: Minimization of Nonlinear Errors}\label{sec: problem}
\begin{figure}[htb]
	\centering\includegraphics[width=0.99\textwidth]{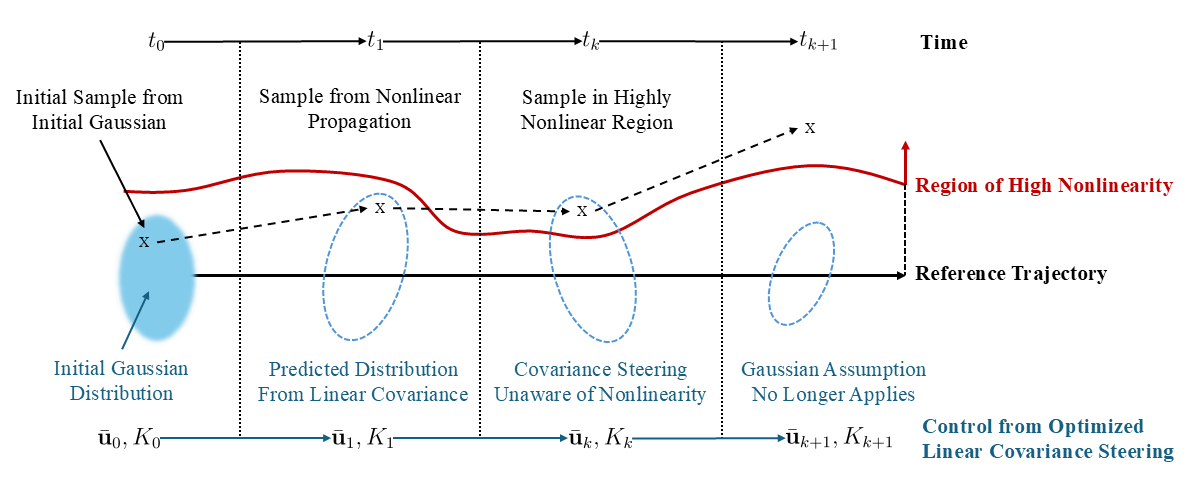}
	\caption{Diagram of Nonlinearity Effects in Linear Covariance Steering.}
	\label{fig: diagram}
\end{figure}

Measures of nonlinearity give insights into how linear approximations hold in nonlinear environments. As mentioned in Section~\ref{sec: intro}, many of the existing MoNs are not formulated for control optimization. For linear covariance steering, understanding how current errors can influence future errors is essential for determining current and forthcoming control actions. Figure~\ref{fig: diagram} illustrates the issues that can arise when nonlinearities are not considered during the control optimization, primarily the divergence from the underlying Gaussian assumptions of the optimization process. Thus, the objectives of this problem are twofold: (1) to determine how nonlinear errors at a specific time impact the trajectory's overall errors, and (2) to ensure that the MoN is structured in a form that can be efficiently solved by optimizers. 

\subsection{Accumulation of Nonlinear Errors}\label{sec: accumulation of NL errors}
Let $\boldsymbol{x}_k^*$ be a reference trajectory where linearization is performed. Then let $\delta\boldsymbol{x}_k$ be the actual state deviation and $\delta\boldsymbol{\tilde{x}}_k$ denote the state deviation propagated through linearized dynamics. At the initial time, there are no nonlinear effects, or $\delta\boldsymbol{x}_0 = \delta\boldsymbol{\tilde{x}}_0$. Let $\delta\boldsymbol{x}^{+}_k$ denote the state deviation after the input of an impulsive control ($\delta\boldsymbol{x}^{+}_k = \delta\boldsymbol{x}_k + B\boldsymbol{u}_k$) and similarly for $\delta\boldsymbol{\tilde{x}}^{+}_k$ ($\delta\boldsymbol{\tilde{x}}^{+}_k = \delta\boldsymbol{\tilde{x}}_k + B\boldsymbol{u}_k$). Currently, the nonlinear contributions from control inputs are neglected. Then, 
\begin{equation}
\begin{aligned}
\delta\boldsymbol{x}_1 = A_0 \delta\boldsymbol{x}^{+}_0 + \mathcal{R}(\delta\boldsymbol{x}^{+}_0) &\qquad&
\delta\boldsymbol{\tilde{x}}_1 = A_0 \delta\boldsymbol{\tilde{x}}^{+}_0 = A_0 \delta\boldsymbol{x}^{+}_0
\end{aligned}
\end{equation}
where  $A_k$ is the linearized dynamics from $t_k$ to $t_{k+1}$ and $\mathcal{R}(\cdot)$ is the ``remainder'' term that is neglected from the linearization.  From Eq.~\eqref{eq: dynamics series expansion}, it can be seen that at any $t_k$, this remainder term is the rest of the series expansion about the reference, or
\begin{equation}\label{eq: linearized remainder}
\mathcal{R}(\delta\boldsymbol{x}^{+}_k) = \sum\limits^{\infty}_{m = 2} \frac{1}{m!} \Phi^{(m)}(t_{k+1},t_k) \cdot (\delta\boldsymbol{x}^{+}_k)^m
\end{equation}

The error associated with linearization at $t_1$ would then be
\begin{equation}
\begin{aligned}
\boldsymbol{\epsilon}_1 = \delta\boldsymbol{x}_1 - \delta\boldsymbol{\tilde{x}}_1 = \mathcal{R}(\delta\boldsymbol{x}^{+}_0)
\end{aligned}
\end{equation}

For $t_2$ and $t_3$,
\begin{equation}
\begin{aligned}
\begin{aligned}
\delta\boldsymbol{x}_2 &= A_1 \delta\boldsymbol{x}^{+}_1 + \mathcal{R}(\delta\boldsymbol{x}^{+}_1) \\
\delta\boldsymbol{\tilde{x}}_2 &= A_1 \delta\boldsymbol{\tilde{x}}^{+}_1 \\
\boldsymbol{\epsilon}_2 &= A_1 \delta\boldsymbol{x}^{+}_1 + \mathcal{R}(\delta\boldsymbol{x}^{+}_1) - A_1 \boldsymbol{\delta\tilde{x}}^{+}_1 \\
           &= \underbrace{A_1 (\delta\boldsymbol{x}^{+}_1 - \delta\boldsymbol{\tilde{x}}^{+}_1)}_{= A_1 (\delta\boldsymbol{x}_1 - \delta\boldsymbol{\tilde{x}}_1) = A_1 \boldsymbol{\epsilon}_1} + \mathcal{R}(\delta\boldsymbol{x}^{+}_1) \\
           &= A_1 \mathcal{R}(\delta\boldsymbol{x}^{+}_0) + \mathcal{R}(\delta\boldsymbol{x}^{+}_1)
\end{aligned}
& \qquad &
\begin{aligned}
\delta\boldsymbol{x}_3 &= A_2 \delta\boldsymbol{x}^{+}_2 + \mathcal{R}(\delta\boldsymbol{x}^{+}_2) \\
\delta\boldsymbol{\tilde{x}}_3 &= A_2 \delta\boldsymbol{\tilde{x}}^{+}_2 \\
\delta\boldsymbol{\epsilon}_3 &= A_2 \delta\boldsymbol{x}^{+}_2 + \mathcal{R}(\delta\boldsymbol{x}^{+}_2) - A_2 \delta\boldsymbol{\tilde{x}}^{+}_2\\
           &= A_2 \boldsymbol{\epsilon}_2 + \mathcal{R}(\delta\boldsymbol{x}^{+}_2) \\
           &= A_2 A_1 \mathcal{R}(\delta\boldsymbol{x}^{+}_0) + A_2 \mathcal{R}(\delta\boldsymbol{x}^{+}_1) + \mathcal{R}(\delta\boldsymbol{x}^{+}_2)
\end{aligned}
\end{aligned}
\end{equation}
and this process can be extended to any $t_k$. It is clear from this analysis that the nonlinear error at $t_k$ ($\boldsymbol{\epsilon}_k$) is a function of the linearized-remainder error of the current timestep ($\mathcal{R}(\cdot)$) as well as previous intervals. Thus, nonlinear errors accumulate throughout the time horizon, so minimizing this error requires knowledge of every nonlinear contribution across the time horizon. Calculating this involves using higher-order state transition tensors. 

\subsection{Bounding Nonlinear Errors and Converting to Convex Form}
As shown in Section~\ref{sec: accumulation of NL errors}, nonlinear errors accumulate over the time horizon, and state transition tensors can be used to quantify these errors. However, tensor operations introduce nonconvexities into the optimization problem, thereby sacrificing some of the favorable convex properties of linear covariance steering. This section introduces a bound to nonlinearities that can transform this problem into a form that is solvable with convex optimization.

Firstly, let's consider minimizing only the $l^2$-norm of the nonlinear error. Using the triangle inequality, each reminder term can be isolated within its own norm.

\begin{equation} \label{eq: NL error bound}
\begin{aligned}
\norm{\boldsymbol{\epsilon}_1}_2 &= \norm{\mathcal{R}(\delta\boldsymbol{x}^{+}_0)}_2 = \norm{I_{n_x} \mathcal{R}(\delta\boldsymbol{x}^{+}_0)}_2 \\
~
\norm{\boldsymbol{\epsilon}_2}_2 &= \norm{A_1 \mathcal{R}(\delta\boldsymbol{x}^{+}_0) + \mathcal{R}(\delta\boldsymbol{x}^{+}_1)}_2  
\leq \norm{A_1 \mathcal{R}(\delta\boldsymbol{x}^{+}_0)}_2 + \norm{\mathcal{R}(\delta\boldsymbol{x}^{+}_1)}_2 \\
~
\norm{\boldsymbol{\epsilon}_3}_2 &= \norm{A_2 A_1 \mathcal{R}(\delta\boldsymbol{x}^{+}_0) + A_2 \mathcal{R}(\delta\boldsymbol{x}^{+}_1) + \mathcal{R}(\delta\boldsymbol{x}^{+}_2)}_2 \leq \norm{A_2 A_1 \mathcal{R}(\delta\boldsymbol{x}^{+}_0)}_2 + \norm{A_2 \mathcal{R}(\delta\boldsymbol{x}^{+}_1)}_2 + \norm{\mathcal{R}(\delta\boldsymbol{x}^{+}_2)}_2 
\end{aligned}
\end{equation}

It can be seen that the right-hand sides of Eq.~\eqref{eq: NL error bound} is a summation of functions of the form $\norm{A \mathcal{R}(\delta\boldsymbol{x}^{+}_k)}_2$ where $A = \prod_i A_i$. Since $\mathcal{R}(\boldsymbol{x}^{+}_k)$ is an infinite series, realistically only STTs up to a certain order are considered. Let $m^* \geq 2$ denote the maximum STT order considered. Then,
\begin{equation}
\mathcal{R}(\delta\boldsymbol{x}^{+}_k) \approx \sum\limits^{m^*}_{m = 2} \frac{1}{m!} \Phi^{(m)}(t_{k+1},t_k) \cdot (\delta\boldsymbol{x}^{+}_k)^m
\end{equation}

Solving for the general case of $\norm{A \mathcal{R}(\delta\boldsymbol{x}^{+}_k)}_2$, the triangle inequality can be used again to separate each $m$-th order STT into its norm.
\begin{equation}\label{eq: pre-bound}
\norm{A \mathcal{R}(\delta\boldsymbol{x}^{+}_k)}_2 \approx \norm{A \sum\limits^{m^*}_{m = 2} \frac{1}{m!} \Phi^{(m)}(t_{k+1},t_k) \cdot (\delta\boldsymbol{x}^{+}_k)^m}_2 \leq 
\sum\limits^{m^*}_{m = 2} \frac{1}{m!} \norm{A \Phi^{(m)}(t_{k+1},t_k) \cdot (\delta\boldsymbol{x}^{+}_k)^m}_2
\end{equation}

The product for $A \Phi^{(m)}(t_{k+1},t_k)$ can be expressed by a tensor $\mathcal{B}^{(m)}$:
\begin{equation}
    \mathcal{B}^{(m)}_{j i_1 i_2 \ldots i_m}(t_{k+1},t_k) = A_{j k} \Phi^{(m)}_{k i_1 i_2 \ldots i_m}(t_{k+1},t_k)
\end{equation}

An upper bound for the expression can be made with the 2-norm of the tensor $\norm{\mathcal{B}}_2$ defined by Kulik et al.~\cite{Kulik-tensorNorm}. 
\begin{equation}
\norm{\mathcal{B}^{(m)}(t_{k+1},t_k) \cdot (\delta\boldsymbol{x}^{+}_k)^m}_2 \leq
\norm{\mathcal{B}^{(m)}(t_{k+1},t_k)}_2 \norm{\delta\boldsymbol{x}^{+}_k}^m_2 
\end{equation}
Numerically, some challenges need to be addressed. This bound involves computing the Z-eigenvector corresponding to the largest eigenvalue of the tensor. \texttt{MATLAB} Tensor Toolbox's shifted symmetric higher-order power method \cite{matlab-eigenvalue} is used to find the Z-eigenpair, which is not guaranteed to converge to the maximum eigenpair. To find the maximum eigenpair, this paper follows the approach suggested by Kulik et al.~\cite{Kulik-tensorNorm}, which uses multiple random initial guesses in the search. While not mathematically rigorous, this strategy yields an eigenpair close to the largest. The primary role of this upper bound is to offer a general insight into the evolution of nonlinear errors, so minor inaccuracies may be acceptable. Note that this bound differs from the previous work by Qi and Oguri~\cite{Qi-Min-NL}, which derives an upper bound using the Cauchy–Schwarz and triangle inequalities. If the user prefers to have a closed-form solution, Qi and Oguri~\cite{Qi-Min-NL}'s method is purely algebraic, and its accuracy does not depend on any Z-eigenpair optimization. However, the bound using $\norm{\mathcal{B}}_2$ is mathematically tighter than that of Qi and Oguri~\cite{Qi-Min-NL} (see Appendix~\ref{appendix: Tensor Upper Bounds}).

This upper bound does not require tensor products, which eliminates the problem's nonconvex elements. For clarity, let's define function $g_k\left(A,m \right) = \norm{\mathcal{B}^{(m)}(t_{k+1},t_k)}_2$. Finally,
\begin{equation}
\norm{A \mathcal{R}(\delta\boldsymbol{x}^{+}_k)}_2 
\leq 
\sum\limits^{m^*}_{m = 2} \frac{1}{m!} \norm{\mathcal{B}^{(m)}(t_{k+1},t_k)}_2 \norm{\delta\boldsymbol{x}^{+}_k}^{m}_2
=
\sum\limits^{m^*}_{m = 2} \frac{1}{m!} g_k(A,m) \norm{\delta\boldsymbol{x}^{+}_k}^{m}_2
\end{equation}

Eq.~\eqref{eq: NL error bound} can then be written as the following:
\begin{equation} \label{eq: upper bound NL error, not Block Matrix}
\begin{aligned}
\norm{\boldsymbol{\epsilon}_1}_2 &\leq 
\sum\limits^{m^*}_{m = 2} \frac{1}{m!} g_0(I_{n_x},m) \norm{\delta\boldsymbol{x}^{+}_0}^{m}_2 = \tilde{\epsilon}_1\\
\norm{\boldsymbol{\epsilon}_2}_2 &\leq
\sum\limits^{m^*}_{m = 2} \frac{1}{m!} g_0(A_1,m) \norm{\delta\boldsymbol{x}^{+}_0}^{m}_2
+ 
\sum\limits^{m^*}_{m = 2} \frac{1}{m!} g_1(I_{n_x},m) \norm{\delta\boldsymbol{x}^{+}_1}^{m}_2   = \tilde{\epsilon}_2\\
\norm{\boldsymbol{\epsilon}_3}_2 &\leq 
\sum\limits^{m^*}_{m = 2} \frac{1}{m!} g_0(A_2 A_1,m) \norm{\delta\boldsymbol{x}^{+}_0}^{m}_2
+ 
\sum\limits^{m^*}_{m = 2} \frac{1}{m!} g_1(A_2,m) \norm{\delta\boldsymbol{x}^{+}_1}^{m}_2
+ 
\sum\limits^{m^*}_{m = 2} \frac{1}{m!} g_2(I_{n_x},m) \norm{\delta\boldsymbol{x}^{+}_2}^{m}_2 = \tilde{\epsilon}_3
\end{aligned}
\end{equation}

Let $\boldsymbol{\Tilde{\epsilon}} = \left[\Tilde{\epsilon}_0 \; \Tilde{\epsilon}_1 \ldots \right]^\top$ be defined as a vector containing the upper bound to $\norm{\boldsymbol{\epsilon}_k}_2 \leq \Tilde{\epsilon}_k$ at each $t_k$. Then, Eq.~\eqref{eq: upper bound NL error, not Block Matrix} can be expressed in a block matrix form for ease of expression.
\begin{equation} \label{eq: upper bound NL error}
\boldsymbol{\Tilde{\epsilon}} = 
\begin{bmatrix}
    \Tilde{\epsilon}_0 \\
    \Tilde{\epsilon}_1 \\ 
    \Tilde{\epsilon}_2 \\ 
    \vdots
\end{bmatrix}
=
\sum\limits^{m^*}_{m = 2} \frac{1}{m!} 
\begin{bmatrix}
    0                   &       0       &      0  &\cdots \\
    g_0\left(I_{n_x},m \right)    &       0       &      0  &\cdots \\
    g_0\left(A_1,m \right) & g_1\left(I_{n_x},m \right) & 0  &\cdots \\
    g_0\left(A_2 A_1,m \right) & g_1\left(A_2,m \right) & g_2\left(I_{n_x},m \right)  &\cdots \\
    \vdots & \vdots  & \vdots  &\ddots \\
\end{bmatrix}
\begin{bmatrix}
    \norm{\delta\boldsymbol{x}^{+}_0}^{m}_2 \\
    \norm{\delta\boldsymbol{x}^{+}_1}^{m}_2 \\ 
    \norm{\delta\boldsymbol{x}^{+}_2}^{m}_2 \\ 
    \vdots
\end{bmatrix}
\end{equation}

The STT is computed alongside the linearization process of the dynamics, meaning $g_k(\cdot)\geq 0$ is a constant and is calculated before the optimization process. The optimization variable $\delta\boldsymbol{x}^{+}_k$ is the only other variable needed to bound the nonlinearities. Since $\norm{\delta\boldsymbol{x}^{+}_k}_2$ is a nonnegative convex function, it can be shown that $\norm{\delta\boldsymbol{x}^{+}_k}^{m}_2$ is also a convex function \cite{Boyd-convex}. This results in $\tilde{\epsilon}_k$ being a summation of convex functions, making it a convex function.

\subsection{Nonlinear Minimization Formulation for Linear Covariance Steering}
It is shown that Eq.~\eqref{eq: upper bound NL error} is a convex function of the optimization variables that can provide an upper bound to the nonlinear dynamical errors. It is important to note that, similar to the linearization and discretization process, the STTs are computed \emph{prior} to the optimization process. While computing STTs is known to be a potentially computationally-taxing process, this process alleviates most of it from the optimizer. Thus, this formulation allows any optimization-based controller to efficiently consider nonlinear effects. 

Typically, the state vector comprises a mixture of different units (e.g., position and velocity), and taking the magnitude of the state vector results in inconsistencies due to the differing units. Thus, it is advantageous to split the contributions from position and velocity separately \cite{Kulik-tensorNorm, Williams-floquet}. This approach allows deviations in position or velocity to be directly mapped to their corresponding nonlinear contributions. Assuming $B = \left[ 0_{3\times 3} \quad I_3\right]^\top$, the $6$-D state becomes $(\delta\boldsymbol{x}_k^{+})^\top = \left[ \delta\boldsymbol{r}^\top_k \quad (\delta\boldsymbol{v}_k^+)^\top \right]^\top$. Then,
\begin{equation} 
\begin{aligned}
\begin{bmatrix}
    \Tilde{\epsilon}_{r,0} \\
    \Tilde{\epsilon}_{r,1} \\ 
    \Tilde{\epsilon}_{r,2} \\ 
    \Tilde{\epsilon}_{r,3}\\ \vdots
\end{bmatrix}
&=
\sum\limits^{m^*}_{m = 2} \frac{1}{m!} 
\begin{bmatrix}
    0                   &       0       &      0  &\cdots \\
    g_{r,0}\left(I_{n_x},m \right)    &       0       &      0  &\cdots \\
    g_{r,0}\left(A_1,m \right) & g_{r,1}\left(I_{n_x},m \right) & 0  &\cdots \\
    g_{r,0}\left(A_2 A_1,m \right) & g_{r,1}\left(A_2,m \right) & g_{r,2}\left(I_{n_x},m \right)  &\cdots \\
    \vdots & \vdots  & \vdots  &\ddots \\
\end{bmatrix}
\begin{bmatrix}
    \norm{\boldsymbol{\delta r}_0}_2^m \\
    \norm{\boldsymbol{\delta r}_1}_2^m \\ 
    \norm{\boldsymbol{\delta r}_2}_2^m \\ \vdots
\end{bmatrix} \\
\end{aligned}
\end{equation}
where $g_{r,k}\left(A,m \right) = \norm{\mathcal{B}_r^{(m)}(t_{k+1},t_k)}_2$. The nonlinear errors for velocity $\Tilde{\boldsymbol{\epsilon}}_{v}$ are computed with $g_{v,k}\left(A,m \right)= \norm{\mathcal{B}_v^{(m)}(t_{k+1},t_k)}_2$ similarly. The tensor $\mathcal{B}^{(m)}$ should exclusively contain either positional or velocity units and be partitioned into submatrices such that
\begin{equation} \label{eq: upper bound NL error, final2}
\begin{rcases}
\left(\mathcal{B}_r^{(m)}\right)_{j i_1 \ldots i_m} &= \left(\mathcal{B}^{(m)}\right)_{j i_1 \ldots i_m} \\
\left(\mathcal{B}_v^{(m)}\right)_{j i_1 \ldots i_m} &= \left(\mathcal{B}^{(m)}\right)_{(j+3) (i_1+3) \ldots (i_m+3)}
\end{rcases}
\quad \text{for} \quad 1 \leq {j,i_1,\ldots,i_m} \leq 3 
\end{equation}

In the context of covariance control, $\norm{\delta\boldsymbol{r}_k}_2$ and $\norm{\delta\boldsymbol{v}^{+}_k}_2$ are stochastic in nature and cannot be treated as deterministic values. Instead, they must be substituted with a statistical quantity that defines the region within which nonlinearities are considered. The upper bound of the quantile from Eq.~\eqref{eq: chance constraints quantile} is used to serve as this metric.
\begin{equation} \label{eq: rtilde and vtilde}
\begin{aligned}
\norm{H_r (\bar{\boldsymbol{x}}^+_k-\boldsymbol{x}_k^*)}_2 + \sqrt{m_{\chi^2}(\epsilon_x,3)}\norm{H_r (P_k^+)^{1/2}}_2 &= \tilde{r}_k\\
\norm{H_v (\bar{\boldsymbol{x}}^+_k-\boldsymbol{x}_k^*)}_2 + \sqrt{m_{\chi^2}(\epsilon_x,3)}\norm{H_v (P_k^+)^{1/2}}_2 &= \tilde{v}_k^+
\end{aligned}
\end{equation}
where $(H_r, H_v)$ extracts either the position or velocity information from the state vector, and $(P_k^+)^{1/2} = \left[(\hat{P}_k^+)^{1/2} \quad \tilde{P}_k^{1/2}\right]$. The MoN objective function can be written in the following form:
\begin{equation} \label{eq: upper bound NL error, final}
\begin{aligned}
\begin{bmatrix}
    \Tilde{\epsilon}_{r,0} \\
    \Tilde{\epsilon}_{r,1} \\ 
    \Tilde{\epsilon}_{r,2} \\ 
    \Tilde{\epsilon}_{r,3}\\ \vdots
\end{bmatrix}
&=
\sum\limits^{m^*}_{m = 2} \frac{1}{m!} 
\begin{bmatrix}
    0                   &       0       &      0  &\cdots \\
    g_{r,0}\left(I_{n_x},m \right)    &       0       &      0  &\cdots \\
    g_{r,0}\left(A_1,m \right) & g_{r,1}\left(I_{n_x},m \right) & 0  &\cdots \\
    g_{r,0}\left(A_2 A_1,m \right) & g_{r,1}\left(A_2,m \right) & g_{r,2}\left(I_{n_x},m \right)  &\cdots \\
    \vdots & \vdots  & \vdots  &\ddots \\
\end{bmatrix}
\begin{bmatrix}
    (\tilde{r}_0)^m \\
    (\tilde{r}_1)^m \\ 
    (\tilde{r}_2)^m \\ 
    \vdots
\end{bmatrix} \\
\end{aligned}
\end{equation}
and a similar procedure can be done with the corresponding velocity quantities of $\Tilde{\epsilon}_{v,k}$ in terms of $\tilde{v}_k^+$ and $g_{v,k}\left(A,m \right)$. The final form of the convex optimization problem is as follows:
\begin{equation}\label{eq: min NL convex problem}
\begin{aligned}
\min_{\bar{\boldsymbol{X}},\hat{\boldsymbol{P}}^{1/2},\bar{\boldsymbol{U}},\boldsymbol{K}} & 
\max_k \left[ (1-\lambda) \tilde{\epsilon}_{r,k} + \lambda \tilde{\epsilon}_{v,k} \right]
&\text{(Minimum Nonlinearity)}
\\
\text{s.t. }
& \text{Eq}.~\eqref{eq: cov steering mean}&\text{(Mean Dynamics)}\\
& \text{Eq}.~\eqref{eq: Phat},\eqref{eq: cov steering cov}&\text{(Covariance Dynamics)}\\
& \text{Eq}.~\eqref{eq: rtilde and vtilde},\eqref{eq: upper bound NL error, final}&\text{(Computing MoN)}\\
\end{aligned}
\end{equation}
where $\max(\cdot)$ denotes the pointwise maximum function of the inputs. The equations for $\boldsymbol{\tilde{\epsilon}}_r,\boldsymbol{\tilde{\epsilon}}_v$ are both convex and a function of state and control, so the pointwise maximum of these convex functions is also a convex function \cite{Boyd-convex}. As a result, it can be applied as a convex optimization problem to solve a control problem that minimizes this measure of nonlinearity. Typically, the objective function for covariance steering is in some form that minimizes the distributional properties of fuel costs \cite{Oguri-Chance-Paper, Naoya-Sequential-Cov-Steering, Jerry-CovSteering, Benedikter-CovSteering}. In the minimum nonlinearity approach, the only modifications are to the objective function, and the rest of the convex linear covariance steering framework remains consistent with prior work.

A tuning parameter $\lambda \in [0,1]$ is introduced to resolve the unit inconsistency when attempting to minimize both $\boldsymbol{\tilde{\epsilon}}_r$ and $\boldsymbol{\tilde{\epsilon}}_v$. A higher $\lambda$ emphasizes minimizing nonlinear effects from control and velocity deviations, whereas a lower $\lambda$ emphasizes those from positional deviations. 
\begin{remark}
The formulation aims to minimize the MoN in the objective, rather than constraining them through convex constraints. This reasoning stems from the MoN’s upper-bounding calculations. Before optimization, it is difficult to establish a meaningful threshold for this measure or to determine whether a particular MoN value leads to improved linear performance. 
\end{remark}
\begin{remark}
A note of caution regarding constraints: while this formulation can handle convex constraints as in the original approach, imposing excessive constraints (e.g., path constraints) may degrade its performance. Constraints are prioritized over minimizing the nonlinearities in the optimizer, which may lead to increased nonlinear errors and discrepancies between the optimized solution and actual behaviors. Ironically, this discrepancy can potentially result in violations of the very path constraints being prioritized in the nonlinear simulation. Similarly, incorporating the additional objective term in Eq.~\eqref{eq: min NL convex problem} results in a multi-objective optimization problem, in which the optimizer no longer exclusively minimizes nonlinearity. Instead, it trades off between competing objectives, potentially yielding solutions that are optimal in the linear sense but perform worse in the true nonlinear system.
\end{remark}
\begin{remark}
It is important to emphasize that this formulation modifies the linear covariance steering framework only through the introduction of a new objective function. Sequential techniques can be applied to update the linearization along the trajectory \cite{Naoya-Sequential-Cov-Steering, Benedikter-CovSteering}, but ultimately, both these approaches and the proposed minimum nonlinearity formulation still rely on the assumption of affine transformations of Gaussian distributions. The formulation in this paper seeks to improve consistency between the linearized and true nonlinear behavior of a covariance steering controller. However, it does not directly address the broader problem of distribution steering in nonlinear systems. In particular, the Gaussian assumption inevitably breaks down under nonlinear dynamics, and a minimum nonlinearity approach alone cannot overcome this fundamental limitation. Qi and Oguri~\cite{Qi-Stat-moment-conference, Qi-Stat-moment} proposed a similar non-Gaussian steering formulation, which overcomes this limitation. Nonetheless, this formulation retains the advantage of being solvable via a computationally efficient convex framework, while also accommodating measurement errors. As a note, this is a capability that Qi and Oguri~\cite{Qi-Stat-moment-conference, Qi-Stat-moment}'s non-Gaussian framework currently lacks.
\end{remark}

\section{Numerical Example: Stationkeeping in Halo Orbits}\label{sec: example}
This paper applies its formulation to a stationkeeping example in a halo orbit. Different halo orbits can possess different degrees of nonlinearity throughout their orbits \cite{Jenson-MeasuresOfNonlinearity, Qi-Min-NL}, and it has been observed that linear covariance controllers can provide inaccurate results in nonlinear systems \cite{Naoya-Sequential-Cov-Steering, Jerry-CovSteering, Qi-Min-NL}. This section compares covariance control on a halo orbit using an objective function that minimizes the MoN presented in this paper in comparison to an objective function found in previous works.

\subsection{Dynamics}
The dynamics are under the circular restricted three-body problem (CR3BP) assumption \cite{Zimovan-Spreen-Halo-Stability}.
\begin{equation}
\begin{aligned}
\ddot{x} - 2\dot{y} = \frac{\partial U}{\partial x}, \qquad \ddot{y} + 2\dot{x} = \frac{\partial U}{\partial y}, \qquad \ddot{z} = \frac{\partial U}{\partial z}
\end{aligned}
\end{equation}
where $x,y,z$ are the nondimensional positions in the synodic frame, $U = \frac{1}{2}\left( x^2 + y^2\right)+\frac{1-\mu}{d}+\frac{\mu}{r}$ is the pseudo-potential function, $\mu$ is the nondimensional mass parameter (Earth-Moon system $\mu \approx 0.0122$), and $d = \sqrt{(x+\mu)^2 + y^2 + z^2}$ and $r = \sqrt{(x-1+\mu)^2 + y^2 + z^2}$. A metric used to quantify the instability of an orbit is with time constant $\tau$ [revs] \cite{Zimovan-Spreen-Halo-Stability}.
\begin{equation}
    \tau\text{ [revs] } = \frac{1}{\text{Re}\left[ \text{Ln}\left( \lambda_{\max}\left[\Phi(t + T,t)  \right] \right) \right]} \frac{1}{T}
\end{equation}
where $\lambda_{\max}\left(\cdot\right)$ returns the largest magnitude eigenvalue of the input, $T$ is the period of the orbit, and $\Phi(t + T,t)$ is the monodromy matrix of the orbit. For stable orbits this value is infinity, and for nearly stable orbits, this value is greater than one. The more unstable an orbit, the lower its associated time constant. This metric helps quantify the stability (or instability) of an orbit for this paper's analysis as more unstable orbits can indicate a greater departure from the linear region around the orbit. 

\subsection{Comparison with a Minimum Covariance Approach}
It is important to show that a minimum nonlinearity approach does not necessarily equate to an approach where the optimized solution is as close to the reference as possible. Previous work by Qi and Oguri~\cite{Qi-Min-NL} has demonstrated that imposing path constraints to remain arbitrarily close to the reference does not necessarily resolve the nonlinearity issue. Another possible point of comparison is between the minimum nonlinearity approach and an objective function based on minimum covariance around the reference.

Jenson and Scheeres~\cite{Jenson-CovMin-traj} observed that minimizing the covariance led to a more Gaussian final distribution. Inspired by Jenson and Scheeres~\cite{Jenson-CovMin-traj}, the form for this paper's minimum covariance analysis minimizes the maximum trace of the positional covariances:
\begin{equation}
    \min_{\bar{\boldsymbol{X}},\hat{\boldsymbol{P}},\bar{\boldsymbol{U}},\boldsymbol{K}} \; \max_k \left(\text{trace}(H_r P_k H_r^\top) \right)
    =
    \min_{\bar{\boldsymbol{X}},\hat{\boldsymbol{P}}^{1/2},\bar{\boldsymbol{U}},\boldsymbol{K}} \; \max_k\left(\; \norm{H_r P_k^{1/2}}_{\text{F}}^2\right)
\end{equation}
where the trace of a matrix can be expressed as a function of the Frobenius norm $||\cdot||_{\text{F}}$ of the decomposition of the matrix \cite{matrix-handbook}. Since $H_r P_k^{1/2}$ is affine, the Frobenius norm is a nonnegative convex function, and the square of a nonnegative convex function is also convex. Thus, the pointwise maximum of these functions is convex and can be used as the objective function for convex optimization \cite{Boyd-convex}. Minimizing covariance supports the notion that as the covariance becomes too large, the distribution can experience nonlinearity to a degree in which the Gaussian assumption is no longer applicable. And thus, minimizing the covariance can lead to a minimization of the covariance's exposure to nonlinear effects. However, the main argument for this comparison is to demonstrate that minimizing covariance is not inherently equivalent to minimizing nonlinearity. 

Firstly, covariance propagation approaches in previous covariance steering formulations \cite{Oguri-Chance-Paper, Naoya-Sequential-Cov-Steering} as well as the work by Jenson and Scheeres~\cite{Jenson-CovMin-traj} use a first-order approximation. Regardless of the value suggested by the first-order approximation, it does not directly reflect the magnitude of higher-order terms. Consequently, a larger covariance does not necessarily imply stronger nonlinearities if the higher-order terms remain small. In such cases, a minimum covariance optimizer may reduce covariance in low-nonlinearity regions while inadvertently increasing covariance in regions of strong nonlinearities, since all covariance is penalized equally regardless of its nonlinear impact. Secondly, these approaches do not account for the accumulation of nonlinear errors. If the distribution were to pass the same region of phase space twice, it is reasonable to expect a smaller covariance on the first pass, since the nonlinear errors introduced then will propagate further downstream than those incurred during the second pass. Nevertheless, this paper's minimum nonlinearity approach still has a strong connection to a minimum covariance approach. Both methods still assume Gaussian distributions and therefore rely solely on the mean and covariance during optimization. The minimum covariance approach treats the deviations at every node equally, whereas the minimum nonlinearity approach assigns appropriate weights that account for nonlinear effects.

\subsection{Simulation Setup}
\begin{figure}[!htb]
  \centering 
    \subcaptionbox
    {Halo Family in Synodic Frame.}
    {\includegraphics[width=0.49\textwidth]{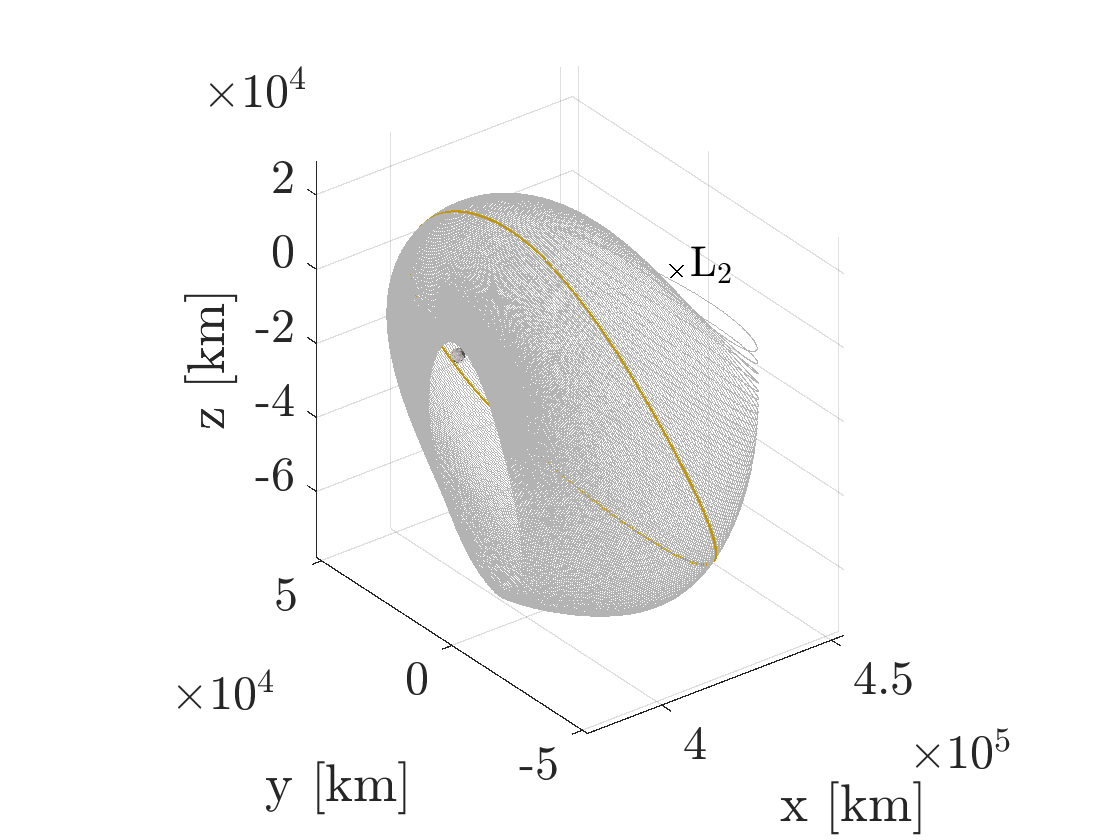}}
    \hskip 0.1truein
    \subcaptionbox
    {Time Constant of Halo Family.}
    {\includegraphics[width=0.49\textwidth]{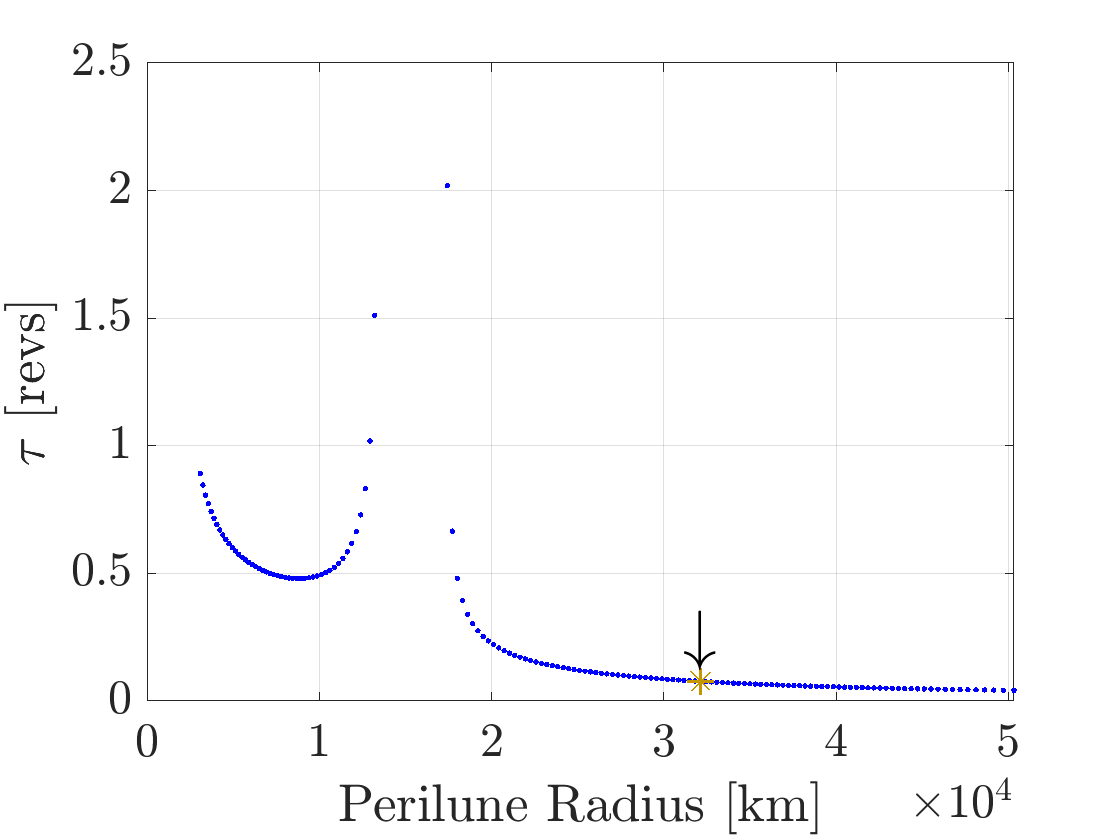}}
    \caption{Southern L$_2$ halo family in barycenter-centered synodic frame. Examined orbit in yellow.}
    \label{fig: halo family and time constant}
\end{figure}

This paper focuses on stationkeeping for the southern Earth-Moon L$_2$ halo family. The orbit chosen is the same orbit as the ``intermediate'' orbit from  Qi and Oguri~\cite{Qi-Min-NL} ($\tau = 0.0770$ revs, $T = 13.071$ days). This is shown in Figure~\ref{fig: halo family and time constant}. 

Table~\ref{tab: CR3BP ICs} lists the parameters for the initial Gaussian distribution at the start of the simulation. Each simulation begins with a mean value $\boldsymbol{\mu}_{0}$ centered at the halo orbit's apolune, and the only constraint during optimization is that the final mean of the linear covariance aligns with this starting value. The halo orbit is used as the reference, meaning $\boldsymbol{x}_k^*$ corresponds to the state along the orbit.

The orbit is discretized linearly in time into nine segments per period, and a $\Delta V$ maneuver occurs only between each segment. The simulation time spans two full revolutions of the orbit, resulting in a total of $N = 19$ nodes. The risk factor $\epsilon_x$ corresponds to a magnitude similar to previous works from Oguri~\cite{Oguri-Chance-Paper}.
\begin{table}[htbp]
	\fontsize{10}{10}\selectfont
    \caption{Parameters for Gaussian Distribution at Initial Time and Constraints}
   \label{tab: CR3BP ICs}
        \centering 
   \begin{tabular}{l  c  c} 
      \hline 
      Parameter & Value & Units\\
      \hline 
      Initial State Mean ($\boldsymbol{\mu}_{0}$) & [1.1300 \quad
         0 \quad
   -0.1767 \quad
         0 \quad
   -0.2255 \quad
         0]$^\top$ & [n.d.] \\
      Initial Position $3\sigma$& 30 &[km]\\
      Initial Velocity $3\sigma$& 3 & [m/s]\\
      Simulation Time & 6.5379 & [n.d.] \\
      Control Magnitude Constraint ($u_{\max}$) & 20 & [m/s] \\
      Risk Factor ($\epsilon_x$) & 0.001 & [n.d.] \\
    \hline
   \end{tabular}
\end{table}

The filter parameters are listed in Table~\ref{tab: filter}. A simple and accurate measurement model is chosen to ensure that the focus remains solely on the controller ($C_k = I_6$, $D_k = \text{blkdiag}(\sigma_{\text{msr},r}I_3, \sigma_{\text{msr},v}I_3)$). Likewise, there are no miscellaneous stochastic perturbations $G_k = 0_{n \times n}$ to keep the dynamics deterministic. These parameters appear in the original formulation from Oguri~\cite{Oguri-Chance-Paper} as well as in modified versions \cite{Naoya-Sequential-Cov-Steering, Jerry-CovSteering} to enhance realism. However, since the primary focus here is on the nonlinearity of a known deterministic dynamical model, the stochastic effects from unmodeled perturbations and filtering are kept to a minimum.
\begin{table}[htbp]
	\fontsize{10}{10}\selectfont
    \caption{Filter Parameters}
   \label{tab: filter}
        \centering 
   \begin{tabular}{l  c  c} 
      \hline 
      Parameter & Value & Units\\
      \hline 
      Initial Position Estimate $3\sigma$& 3 &[km]\\
      Initial Velocity Estimate $3\sigma$& 3 & [m/s]\\
      Position Measurement Uncertainty $(\sigma_{\text{msr},r})$ & 1 & [m] \\
      Velocity Measurement Uncertainty $(\sigma_{\text{msr},v})$ & 10 & [cm/s] \\
    \hline
   \end{tabular}
\end{table}

The tuning parameters of the minimum nonlinearity formulation are defined as follows. Cases considering second- and third-order STTs, corresponding to maximum orders $m^* = 2$ and $m^* = 3$ are analyzed. The weighting parameter $\lambda=0.52$ is selected to be approximately one-half, since the MoNs for position and velocity are nondimensionalized and of equal magnitude in CR3BP. The STT and tensor norm computations took around 30 seconds for $m^* = 2$ and 2 minutes for $m^* = 3$ with \texttt{MATLAB} R2024b. The convex problems are solved via \texttt{CVX} with \texttt{MOSEK} and took around 10 seconds to solve for both cases.

\subsection{Effects on Predicted Distribution}
\begin{figure}[htb]
	\centering\includegraphics[width=0.9\textwidth]{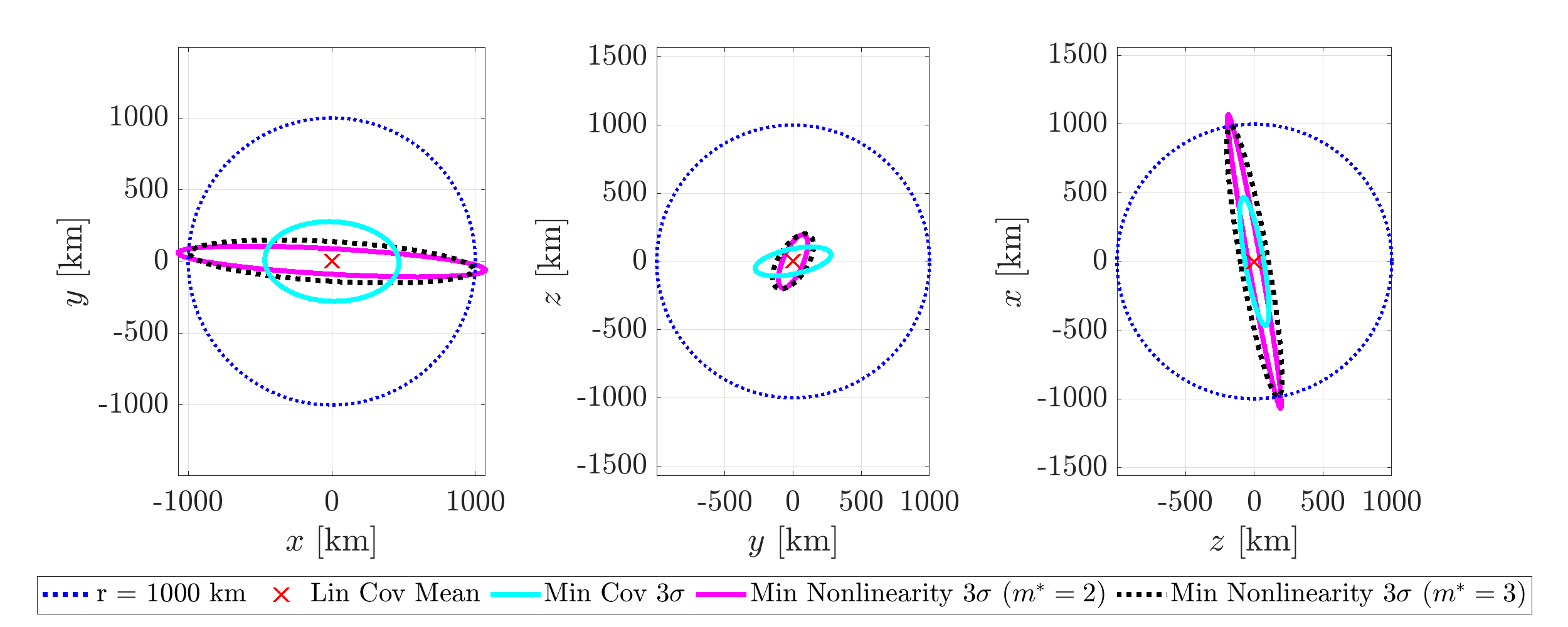}
	\caption{Comparison of the final linear covariance between the minimum covariance and minimum nonlinearity approaches. \emph{Origin normalize to halo orbit.}}
	\label{fig: lin cov compare}
\end{figure}

\begin{figure}[!htbp]
  \centering 
    \subcaptionbox
    {With Minimum Covariance.\label{fig: Results, minCov}}
    {\includegraphics[width=0.9\textwidth]{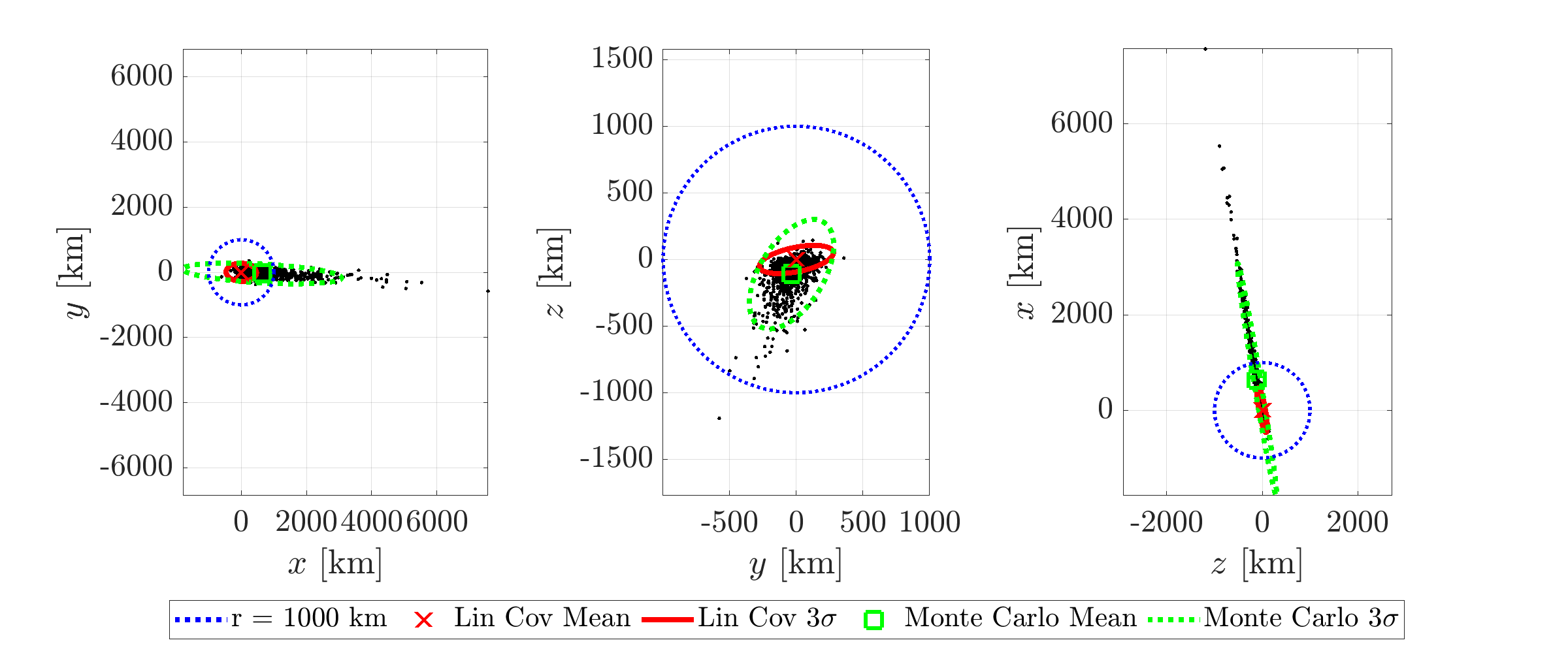}}
    \hskip 0.1truein
    ~
    \subcaptionbox
    {With Minimum Nonlinearity ($m^*=2$).\label{fig: Results, minNL}}
    {\includegraphics[width=0.9\textwidth]{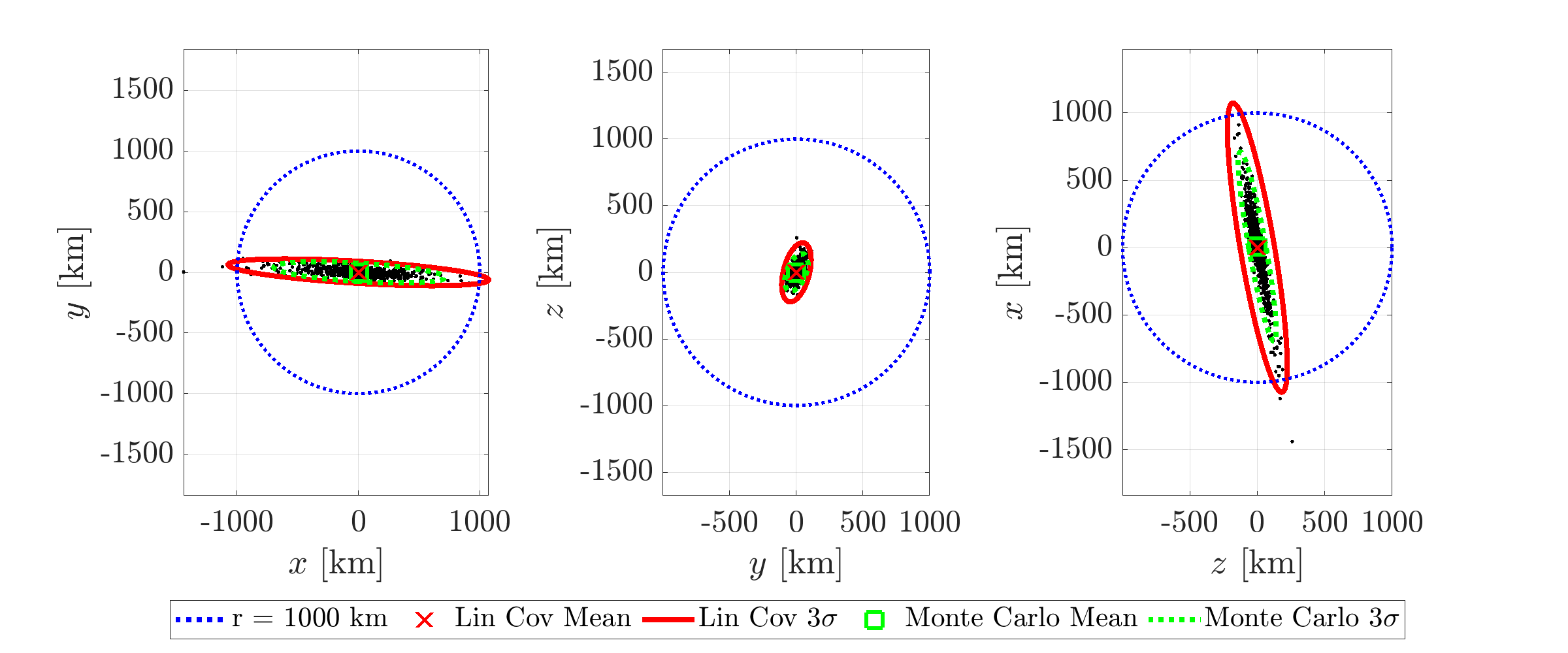}}
    ~
    \subcaptionbox
    {With Minimum Nonlinearity ($m^*=3$).\label{fig: Results, minNL_m3}}
    {\includegraphics[width=0.9\textwidth]{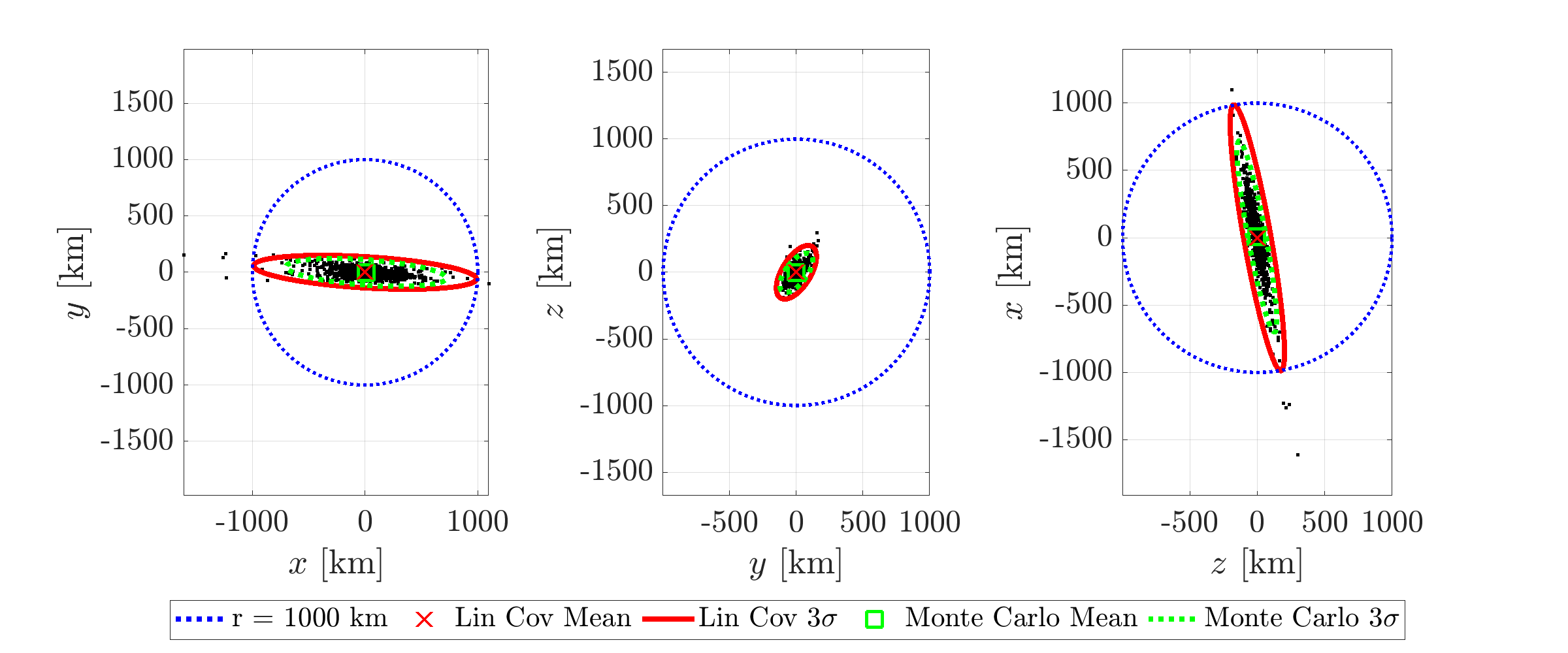}}
    \caption{Monte Carlo ($n_{\text{samples}}=1,000$): Distribution at the terminal time comparing the minimum covariance and minimum nonlinearity approaches. \emph{Origin normalize to halo orbit.}}
    \label{fig: Results}
\end{figure}

Figure~\ref{fig: lin cov compare} shows the predicted linear solution to the optimization problem. The dashed blue line that marks a radius of $r = 1,000$ km serves as a scale for reference when comparing the distribution in this and future figures. Since the mean is constrained, the optimized mean is the same between the minimum covariance and minimum nonlinearity approaches. The difference lies in the optimized linear covariance. Firstly, the minimum covariance inherently seeks to reduce the maximum size of the covariance, and as shown, it has a smaller $3\sigma$ than that of the minimum nonlinearity cases. Secondly, there is a strong similarity between the covariances of $m^*=2$ and $m^*=3$. This suggests that increasing the tensor order does not significantly impact the characterization of nonlinearity in this problem. This observation is consistent with prior work showing that tensor-based MoNs can exhibit similar behavior across different orders for halo orbits \cite{Jenson-MeasuresOfNonlinearity}, but this may not generalize to all problems. In some cases, different orders of STTs can express different trends \cite{calkins-STT}, and the appropriate choice is left to the discretion of the user.

Nevertheless, the behavior predicted by the linear covariance does not guarantee a similar nonlinear behavior. Figure~\ref{fig: Results} presents the effects of the differing control policies in a nonlinear simulation. It is clear from Figure~\ref{fig: Results, minCov} that the linear covariance prediction is no longer representative of the distribution from the Monte Carlo. Both the mean and the covariance no longer match the values predicted by the linear analysis, and in fact appear to diverge from the predicted solution. On the other hand, Figure~\ref{fig: Results, minNL} demonstrates significantly improved agreement between the statistics predicted by the linear model and those obtained from the Monte Carlo distribution. Furthermore, the dispersion of the Monte Carlo distribution is notably smaller than that observed under the minimum covariance approach. Figure~\ref{fig: Results, minNL_m3} shows a similar trend to the lower-order tensor case, where the nonlinear behavior is accurately captured by the linear prediction.

\begin{figure}[!htbp]
  \centering 
    \subcaptionbox
    {With Minimum Covariance.\label{fig: quantile, minCov}}
    {\includegraphics[width=0.48\textwidth]{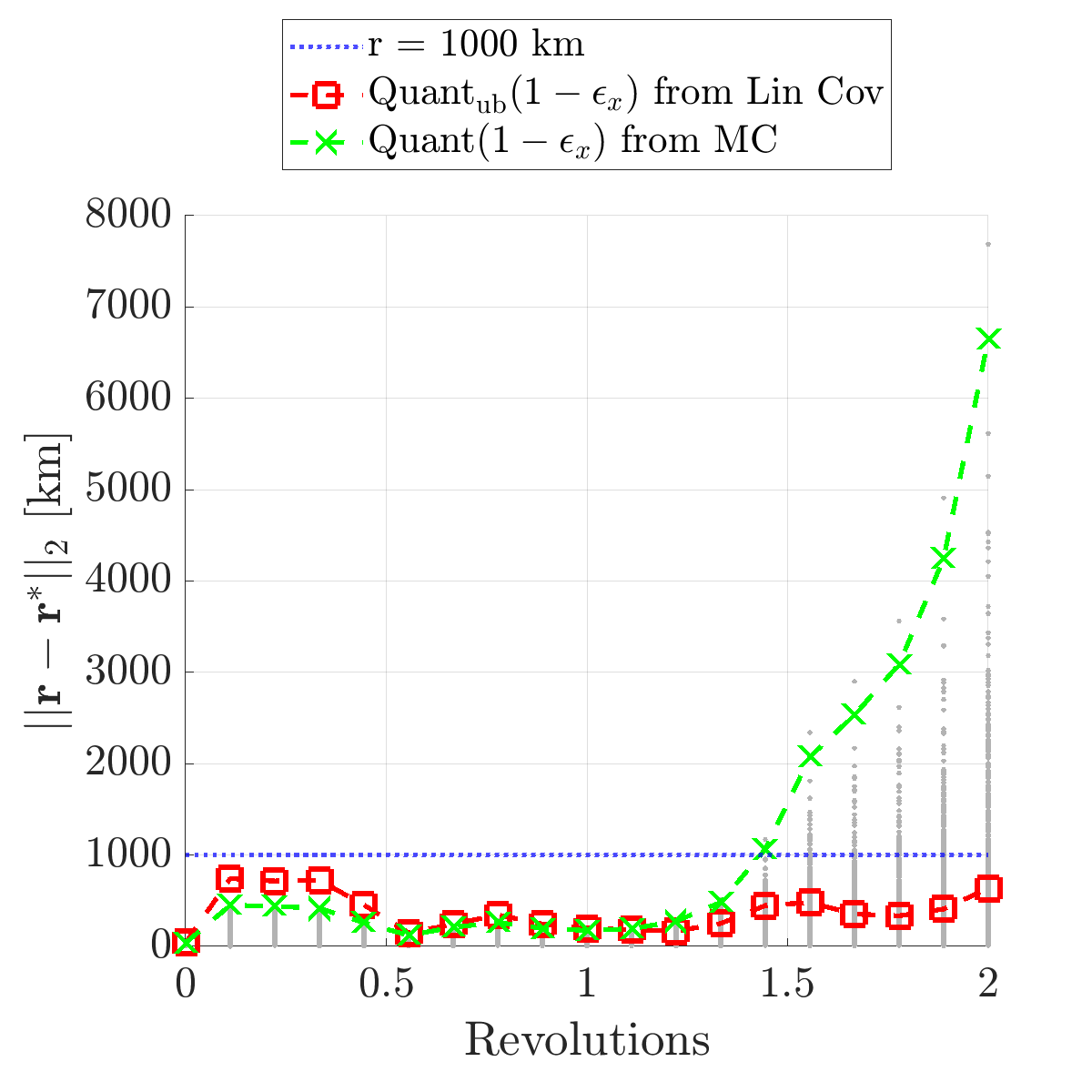}}
    \hskip 0.1truein
    ~
    \subcaptionbox
    {With Minimum Nonlinearity ($m^*=2$).\label{fig: quantile, minNL}}
    {\includegraphics[width=0.48\textwidth]{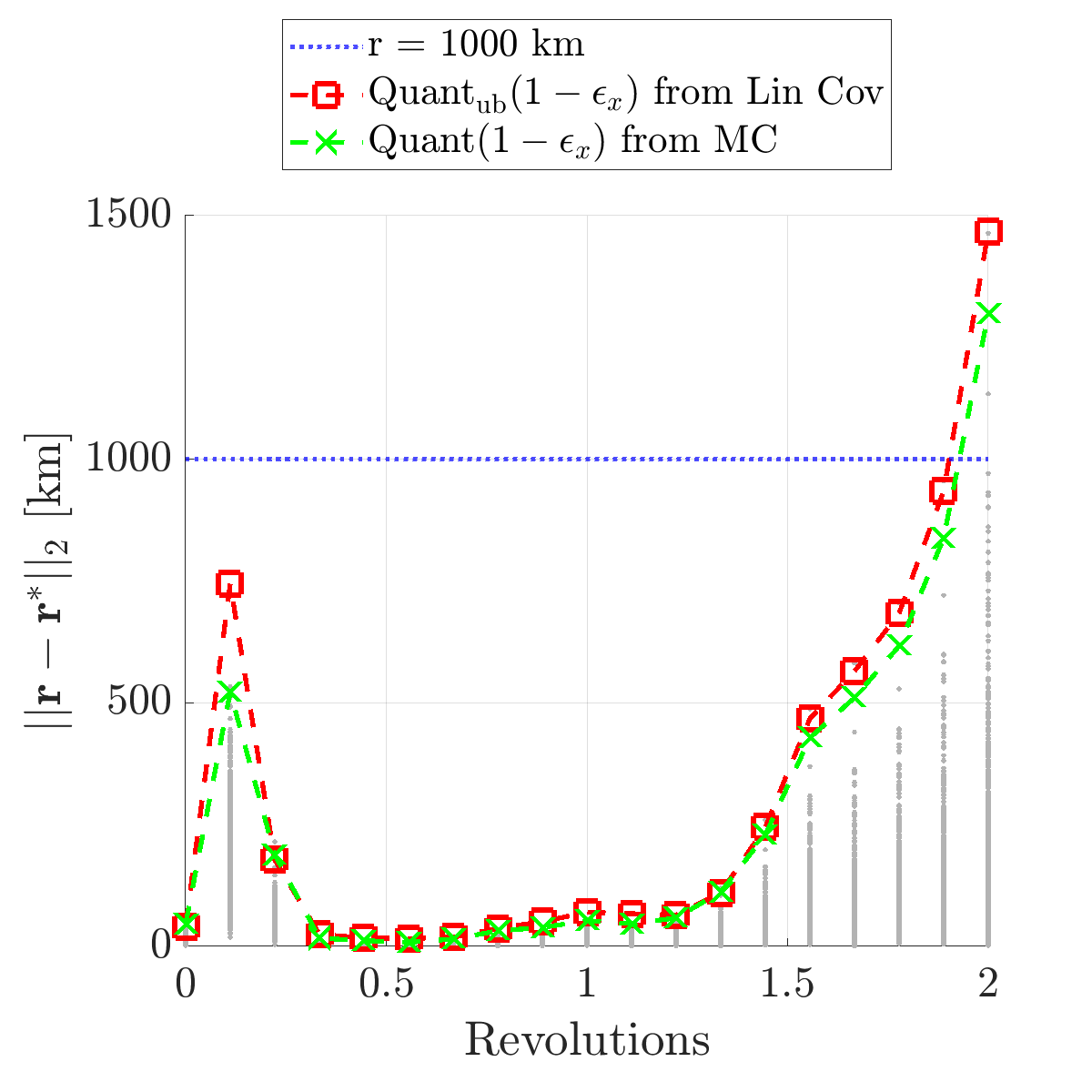}}
    \caption{Monte Carlo ($n_{\text{samples}}=1,000$): Comparison of the upper-bound quantile function $\text{Quant}_{\text{ub}}(\cdot)$ for state deviations predicted by the linear covariance with the actual quantile function $\text{Quant}(\cdot)$ obtained from Monte Carlo simulations.}
    \label{fig: quantile}
\end{figure}

Figure~\ref{fig: quantile} shows the quantile information of state deviations from the reference between the two methods. Each plot contains the upper bound of the positional quantile (calculated by $\tilde{r}_k$ in Eq.~\eqref{eq: rtilde and vtilde}, but denoted as $\text{Quant}_{\text{ub}}$ for ease of reference) compared to the quantile from the Monte Carlo (denoted by $\text{Quant}$). It should be remembered that $\text{Quant}_{\text{ub}}$ is only an upper bound under the Gaussian assumption, and the mean and covariance are computed using linear covariance analysis.

An interesting observation is the linear covariance solution of each method. Regions near perilune for halo orbits generally have stronger nonlinearities than regions around apolune \cite{Jenson-MeasuresOfNonlinearity}. Compare Figure~\ref{fig: quantile, minCov} and \ref{fig: quantile, minNL}, the minimum nonlinearity approach is seen to reduce the covariance rapidly around the perilune region (0.5 revs). While the minimum covariance approach achieves a similar reduction, its reduction is more gradual. 

As for the comparison with the nonlinear simulations, the predicted and actual quantile aligns well at the beginning for both cases. This is likely since there are low nonlinear errors and the distributions are generally Gaussian. For the minimum covariance approach in Figure~\ref{fig: quantile, minCov}, discrepancies start to arise around one revolution in the orbit, which only compounds over time. In contrast, for the minimum nonlinearity approach with $m^*=2$ shown in Figure~\ref{fig: quantile, minNL}, the upper bound accurately captures the Monte Carlo quantile throughout the simulation. The corresponding plot for the minimum nonlinearity approach with $m^*=3$ exhibits a similar trend, but is omitted for brevity. Strictly speaking, the actual distribution for both approaches becomes non-Gaussian due to the nonlinear transformations, but the one from the minimum nonlinearity formulation appears much closer to the predicted Gaussian distribution. This indicates that, for mission operations, the linear covariance derived from the minimum nonlinearity approach provides a more reliable approximation when applied to nonlinear regimes.

\begin{figure}[htb]
	\centering\includegraphics[width=0.6\textwidth]{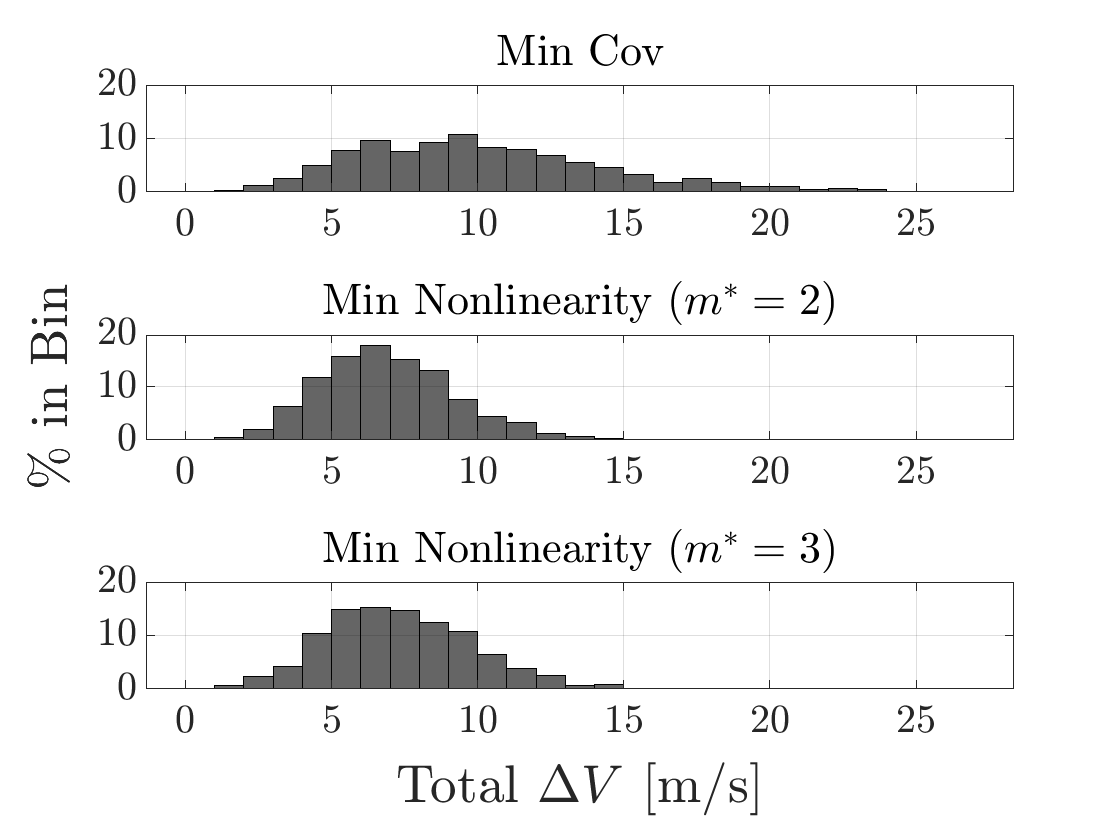}
	\caption{Monte Carlo ($n_{\text{samples}}=1,000$): Histogram of total $\Delta V$ between the methods.}
	\label{fig: deltaV}
\end{figure}

A summary of the $\Delta V$ costs for each method is provided in Figure~\ref{fig: deltaV}. Due to the similarities between $m^*=2$ and $m^*=3$, their corresponding fuel cost distributions are also similar. While the minimum covariance approach yields costs of the same order of magnitude, it exhibits significantly greater dispersion. This is likely due to the increased state dispersion observed in the nonlinear simulations, as shown in Figure~\ref{fig: Results}, which results from poorer control of the nonlinear distribution. Consequently, additional control effort is required, as the optimized gain is less well-suited to the true nonlinear dynamics.

\section{Conclusion}\label{sec: conclusion}
This paper presents a novel formulation that increases the performance of linear covariance controllers in nonlinear environments through the minimization of nonlinear errors. Firstly, it is shown that dynamical errors accumulate with each control action due to the linearization of dynamics. This highlights that considering immediate nonlinearities from previous measures of nonlinearities may not be sufficient, since early errors can propagate and lead to significant downstream inaccuracies. Next, the nonlinear errors are quantified by upper-bounding the contributions using high-order state transition tensors. This paper derives an upper bound as a function of state and control input, and this formulation is notably in a convex form. As a result, it can be used in convex optimization to minimize the nonlinear dynamical errors of linear covariance controllers. Lastly, this paper applies the formulation to a halo orbit stationkeeping example. The results demonstrate an improvement in capturing the nonlinear distributions with linear covariance compared to the existing methods for handling system nonlinearities. This indicates that by using the formulation presented in this paper, the linear covariance controller operates more effectively as intended in nonlinear environments.

\section{Appendix}
\subsection{Comparison of Upper Bounds}\label{appendix: Tensor Upper Bounds}
Given a tensor $\mathcal{T}^{(m+1)}$ with $m+1$ number of modes and dimension $n$ along each mode, and a vector $\boldsymbol{v}$ of dimension $n$, the goal is to get an accurate convex upper bound of the quantity:
\begin{equation}
    \norm{\mathcal{T}^{(m+1)} \cdot \boldsymbol{v}^m}_2 
\end{equation}

Previous works by Qi and Oguri~\cite{Qi-Min-NL} derive a closed-form upper bound using the Cauchy–Schwarz and triangle inequalities. This bound always holds algebraically:
\begin{equation}\label{eq: appendix ub, Qi}
\norm{\mathcal{T}^{(m+1)} \cdot \boldsymbol{v}^m}_2 
\leq 
n^{\frac{m}{2}} \norm{\boldsymbol{v}}^m_2 \norm{\boldsymbol{\mathcal{T}}_{1:\max}}_2
\end{equation}
where $\boldsymbol{\mathcal{T}}_{1:\max}$ is a vector containing the maximum absolute value of $\mathcal{T}^{(m+1)}$ throughout its first mode. This is equivalent to the $j$-th element of $\boldsymbol{\mathcal{T}}_{1:\max}$ expressed as 
\begin{equation}
\left( \boldsymbol{\mathcal{T}}_{1:\max} \right)_j = \max_{i_2 \ldots i_{m+1}} \left( \; \abs{\mathcal{T}_{j i_2 \ldots i_{m+1}}} \;\right)
\end{equation}

Recently, Kulik et al.~\cite{Kulik-tensorNorm} propose a 2-norm for an arbitrary tensor in which an upper bound can be constructed:
\begin{equation}\label{eq: appendix ub, Kulik}
\norm{\mathcal{T}^{(m+1)} \cdot \boldsymbol{v}^m}_2 
\leq 
\norm{\mathcal{T}^{(m+1)}}_2 \norm{\boldsymbol{v}}^m_2
\end{equation}
As discussed earlier, tensor eigenpair algorithms are not guaranteed to converge to the largest eigenvalue. However, they often do so in practice, and so Kulik et al.~\cite{Kulik-tensorNorm} recommend using multiple random initializations when searching for the largest eigenvalue.

A numerical demonstration is made to compare the tightest of each bound. With $n =3$ and $m=2$, each element in $\mathcal{T}^{(m+1)}_{\text{sample}}$ and $\boldsymbol{v}_{\text{sample}}$ is sampled from a uniform distribution from $[-100,100]$. A percent difference is used to determine the tightness of each bound:
\begin{equation}
    \text{$\%$ difference} = \frac{\text{(upper bound)}-
    \norm{\mathcal{T}^{(m+1)}_{\text{sample}} \cdot \boldsymbol{v}^m_{\text{sample}}}_2}{\norm{\mathcal{T}^{(m+1)}_{\text{sample}} \cdot \boldsymbol{v}^m_{\text{sample}}}_2} \times 100
\end{equation}
where the upper bound corresponds to either Qi's or Kulik's method in Eq.~\eqref{eq: appendix ub, Qi} or Eq.~\eqref{eq: appendix ub, Kulik} respectively. Ten random initial guesses with \texttt{eig\_sshopm()} from \texttt{MATLAB} Tensor Toolbox \cite{matlab-eigenvalue} are used for Kulik's calculation of the tensor 2-norm.  

\begin{figure}[!htb]
	\centering\includegraphics[width=0.6\textwidth]{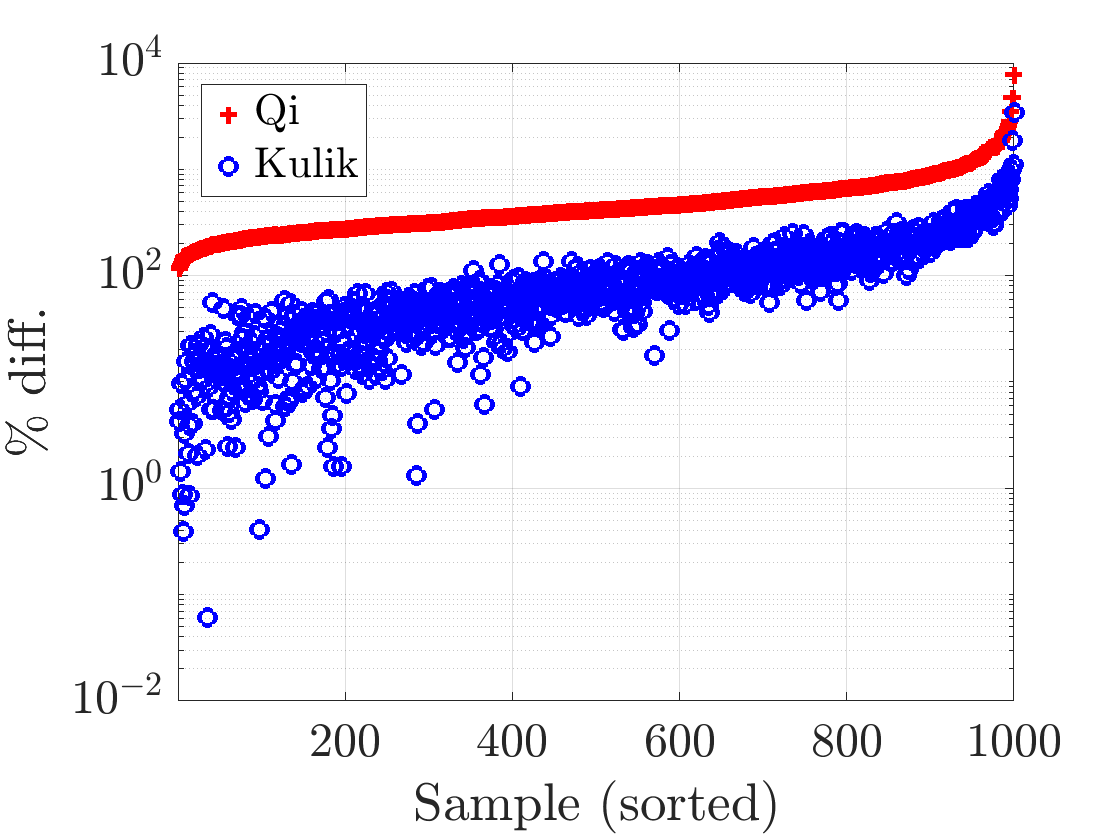}
	\caption{Monte Carlo ($n_{\text{samples}}=1,000$): Comparison of Inequality Tightness. \emph{Results sorted by Qi's upper bound.}}
	\label{fig: Bnorm Compare}
\end{figure}
From Figure~\ref{fig: Bnorm Compare}, it can be seen that Kulik's upper bound is nearly two orders of magnitude tighter than Qi's upper bound in every case. Note that if the percent difference is negative, the upper bound is violated. Although Kulik's upper bound can, in principle, be violated if the maximum tensor eigenpair is not found, no such violations were observed in this demonstration.

\section*{Acknowledgments}
The authors thank Dr. Jackson Kulik for his recommendations on the bound for Eq.~\eqref{eq: pre-bound} in the form of a tensor 2-norm. A large language model was used solely for language editing to improve concision and readability. It was not used for generating content, analysis, or conclusions. All outputs were reviewed and validated by the authors. This material is based upon work supported by the Air Force Office of Scientific Research under award number FA9550-23-1-0512. 

\bibliography{references}

\end{document}